\documentclass[a4paper,12pt]{article}
\usepackage{amscd}
\usepackage{amssymb}                                    %
\usepackage{mathrsfs}                    %
\usepackage{amsmath}                    %
\usepackage{amsfonts}                   %
\usepackage{amsthm}                     %
\usepackage[latin1]{inputenc}           %
\usepackage[arrow, matrix, curve,cmtip]{xy}    %
\usepackage[english]{babel}               
\usepackage{graphics}
\usepackage[unicode,hypertexnames]{hyperref}
\usepackage{bbm}
\hypersetup{
    urlcolor = red,
	colorlinks = true,
	linkcolor = blue,
	citecolor = blue,
    linktocpage = true,
	pdftitle = {A Smooth Model for the String Group},
    pdfauthor = {Thomas Nikolaus, Christoph Sachse and Christoph Wockel},
    bookmarksopen = true,
    bookmarksopenlevel = 1,
	unicode = true,
    hypertexnames = false
  }

\theoremstyle{plain}
\newtheorem{theorem}{Theorem}[section]
\newtheorem*{nntheorem}{Theorem}
\newtheorem{corollary}[theorem]{Corollary}
\newtheorem{proposition}[theorem]{Proposition}
\newtheorem{lemma}[theorem]{Lemma}
\theoremstyle{definition}
\newtheorem{example}[theorem]{Example}
\newtheorem{remark}[theorem]{Remark}
\newtheorem{problem}[theorem]{Problem}
\newtheorem{definition}[theorem]{Definition}

\newcommand{\cat}[1]{\ensuremath{\boldsymbol{\op{#1}}}}

\newcommand{\op}[1]{\ensuremath{\operatorname{#1}}}
\newcommand{\ol}[1]{\ensuremath{\overline{#1}}}
\newenvironment{frcseries}{\fontfamily{frc}\selectfont}{}

\DeclareMathAlphabet{\mathpzc}{OT1}{pzc}{m}{it}
\DeclareMathAlphabet{\mathcalligra}{T1}{calligra}{m}{n}

\DeclareFontFamily{OT1}{pzc}{}
\DeclareFontShape{OT1}{pzc}{m}{it}%
             {<-> s * [1.150] pzcmi7t}{}
\DeclareMathAlphabet{\mathscr}{OT1}{pzc}%
                                 {m}{it}
\newcommand{\Gau}{\ensuremath{\mathscr{Gau}}}
\newcommand{\Gaut}{\ensuremath{\widehat{\mathscr{Gau}}}}
\newcommand{\Gautt}{\ensuremath{\widetilde{\mathscr{Gau}}}}

\newcommand{\Gauc}{\ensuremath{\mathscr{Gau}^c}}
\newcommand{\gau}{\ensuremath{\mathfrak{gau}}}
\newcommand{\Aut}{\mathrm{Aut}}
\newcommand{\Homeo}{\mathrm{Homeo}}
\newcommand{\String}{\ensuremath{\mathscr{String}_{G}}}
\newcommand{\BCString}{\ensuremath{\mathcal{P}\mathscr{String}_{G}}}
\newcommand{\Spring}{\ensuremath{\mathscr{String}}}
\newcommand{\Stringc}{\ensuremath{\mathscr{String}_{G}^c}}
\newcommand{\Diff}{\mathrm{Diff}}
\newcommand{\fk}{\mathfrak{k}}
\newcommand{\fg}{\mathfrak{g}}
\newcommand{\Z}{\mathbbm{Z}}
\newcommand{\R}{\mathbbm{R}}
\newcommand{\calk}{\mathcal{K}}
\newcommand{\calh}{\mathcal{H}}
\newcommand{\calg}{\mathcal{G}}
\newcommand{\calb}{\mathcal{B}}

\newcommand{\STRING}{\ensuremath{\op{\textnormal{\texttt{STRING}}}_{G}}}
\newenvironment{tabsection}{}{}
\newcommand\iso{\xrightarrow{\sim}}
\newcommand{\from}{\ensuremath{\nobreak\colon\nobreak}}

\newcommand{\upi}{\ensuremath{\underline{\pi}}}
\newcommand{\Hom}{\mathrm{Hom}}
\newcommand{\id}{\mathrm{id}}

\def\zwoabun#1#2{2\text{\text{-}}\buntech{#1}{}(#2)}
\def\buntech#1#2{\mathcal{B}\hspace{-0.01em}un^{#2}_{\hspace{-0.04em}#1}}
\def\bun#1#2{\buntech{#1}{}(#2)}

\begin{document}
\title{A Smooth Model for the String Group}
\author{
Thomas Nikolaus\\\small{}Universit\"at Regensburg, Fakult\"at für Mathematik\\\small{}Universitätsstra\ss{}e 31, D-93053 Regensburg\\\small{}Thomas1.Nikolaus@mathematik.uni-regensburg.de 
\and
Christoph Sachse\\\small{}Universit\"at Hamburg, Fachbereich Mathematik\\\small{}Bundesstra\ss{}e 55, D-20146 Hamburg\\\small{}chsachse@googlemail.com
\and
Christoph Wockel\\\small{}Universit\"at Hamburg, Fachbereich Mathematik\\\small{}Bundesstra\ss{}e 55, D-20146 Hamburg\\\small{}christoph.wockel@math.uni-hamburg.de
}
\date{}
\maketitle

\begin{abstract}
 \noindent We construct a model for the string group as an
 infinite-dimensional Lie group. In a second step we extend this model
 by a contractible Lie group to a Lie 2-group model. To this end we need
 to establish some facts on the homotopy theory of Lie 2-groups.
 Moreover, we provide an explicit comparison of string structures for
 the two models and a uniqueness result for Lie 2-group models.
\end{abstract}

\tableofcontents

\section{Introduction}

\begin{tabsection}
 String structures and the string group play an important role in algebraic
 topology \cite{hoplect, LurSurvey, BunNau09}, string theory \cite{kill87,
 FreedMoore} and geometry \cite{wit88, Stolz}. The group $\Spring$ is defined
 to be a 3-connected cover of the spin group or. More generally, we denote by
 $\String$ the 3-connected cover of any compact, simple and 1-connected
 Lie group $G$ \cite{StolzTeichner04What-is-an-elliptic-object}. This
 definition fixes only its homotopy type and makes abstract homotopy theoretic
 constructions possible. These models are not very well suited for geometric
 applications, one is rather interested in concrete models that carry, for
 instance, topological or even Lie-group structures.
 
 There is a direct cohomological argument showing that $\String$ cannot be a
 finite $CW$-complex or a finite-dimensional manifold (see Corollary
 \ref{cor:string-is-not-finite-dimensional}), so the best thing one can hope
 for is a model for $\String$ as a topological group or an infinite-dimensional
 Lie group. There have been various constructions of models of $\String$ as
 $A_{\infty}$-spaces or topological groups, but the question whether an
 infinite-dimensional Lie group model is also possible remains open. One of the
 main contributions of the present paper is to give an affirmative answer to
 this question and provide an explicit Lie group model, based on a topological
 construction of Stolz \cite{Stolz}.
 
 Something that is not directly apparent from the setting of the problem is
 that string group models as Lie 2-groups are something more natural to expect
 when taking the perspective of string theory or higher homotopy theory into
 account. However, the notion of a Lie 2-group model deserves a thorough
 clarification itself. We discuss this notion carefully by establishing the
 relevant homotopy theoretic facts about infinite-dimensional Lie 2-groups and
 promote our Lie group model $\String$ to such a Lie 2-group model $\STRING$.
 \\
 
 Before we outline our construction let us briefly summarize the existing ones.
 Let $G$ be throughout a compact, simple and 1-connected Lie group. One
 model for $\String$ can be obtained from pulling back the path fibration
 $P K(\Z,3)\to K(\Z,3)$ along a characteristic map $u\from G\to K(\Z,3)$. This
 is a standard construction of the Whitehead tower and leads to a model of
 $\String$ as a space. Since this construction also works for a characteristic
 map $BG\to K(\Z,4)$, each 3-connected cover is homotopy equivalent to a loop
 space and thus admits an $A_{\infty}$-structure. Taking a functorial and
 product preserving construction of the Whitehead tower one even obtains a
 model as a topological group. Unfortunately, these models are not very
 tractable.
 
 There are more geometric constructions of $\String$, for instance the one by
 Stolz in \cite{Stolz}. The model given there has as an input the basic
 principal $PU(\calh)$-bundle $P$ over $G$, where $\cal{H}$ is a complex
 separable infinite-dimensional Hilbert space. Stolz then defines a model for
 $\String$ as a topological group together with a homomorphism $\String\to G$
 whose kernel is the group of continuous gauge transformations of the bundle
 $P$. Our constructions will be based on this idea. In
 \cite{StolzTeichner04What-is-an-elliptic-object} Stolz and Teichner construct
 a model for $\String$ as an extension of $G$ by $PU(\calh)$. It is a natural
 idea to equip this model with a smooth structure. But this does not work since
 this extension is constructed as a pushout along a positive energy
 representation of the loop group of $G$ which is not smooth.
 
 We now come to Lie 2-group models. One such construction has been given by
 Henriques \cite{Henriques08Integrating-Lsb-infty-algebras}, based on work of
 Getzler \cite{Getzler09Lie-theory-for-nilpotent-Linfty-algebras}. The basic
 idea underlying this construction is to apply a general integration procedure
 for $L_\infty$-algebras to the string Lie 2-algebra. To make this construction
 work one has to weaken the naive notion of a Lie 2-group and moreover work in
 the category of Banach spaces. Similarly, the model of Schommer-Pries
 \cite{Schommer-Pries10Central-Extensions-of-Smooth-2-Groups-and-a-Finite-Dimensional-String-2-Group}
 realizes $\String$ as a stacky Lie 2-group, but it has the advantage of being
 finite-dimensional. This model is constructed from a cocycle in Segal's
 cohomology of $G$ \cite{Segal70Cohomology-of-topological-groups}.
 
 A common weakness of the above Lie 2-group models is that they are not strict,
 i.e. not associative on the nose but only up to an additional coherence. This
 complication is not present in the strict 2-group model of Baez, Crans,
 Schreiber and Stevenson from
 \cite{BaezCransStevensonSchreiber07From-loop-groups-to-2-groups}. It is
 constructed from a crossed module $\widehat{\Omega G}\to P_{e}G$, built out of
 the level one Kac-Moody central extension $\widehat{\Omega G}$ of the based
 smooth loop group $\Omega G$ of $G$ and its path space $P_{e}G$. The price to
 pay is that the model is infinite dimensional, but the strictness makes the
 corresponding bundle theory more tractable \cite{NikolausWaldorf11}.\\
 
 Summarizing, quite some effort has been made in constructing models for
 $\String$ that are as close as possible to finite-dimensional Lie groups.
 However, one of the most natural questions, namely whether there exists an
 infinite-dimensional \emph{Lie group} model for $\String$ is still open. We
 answer this question by the following result.
 
 Let $P\to G$ be a basic smooth principal $PU(\calh)$-bundle, i.e.,
 $[P]\in [G,BPU(\calh)]\cong H^{3}(G,\Z) = \Z$ is a generator. In Section
 \ref{sect:gaugegrp_lie} we review the fact that $\Gau(P)$ is a Lie group
 modeled on the infinite-dimensional space of vertical vector fields on $P$.
 The main result of Section \ref{sec:group} is then the following enhancement
 of the model from \cite{Stolz} to a Lie group model.
 \begin{nntheorem}[Theorem
  \ref{mainthm1}] Let $G$ be a compact, simple and 1-connected Lie group,
  then there exists a model $\String$ of the string group, which is a
  metrizable Fr\'echet--Lie group and turns
  \begin{equation*}
   \Gau(P)\to \String \to G
  \end{equation*}
  into an extension of Lie groups. It is uniquely determined up to isomorphism
  by this property.
 \end{nntheorem}
 From now on $\String$ will always refer to this particular model.
 Metrizability makes the homotopy theory that we use in the sequel work due to
 results of Palais \cite{pal66}.
 
 In Section \ref{sect:2groups} we introduce the concept of Lie 2-group models
 culminating in Definition \ref{def:twomodel}. An important construction in
 this context is the geometric realization that produces topological groups
 from Lie 2-groups. We show that geometric realization is well-behaved under
 mild technical conditions, such as metrizability.
 
 In Section \ref{sect:the_string_group_as_a_2_group} we then construct a
 central extension $U(1)\to \Gaut(P)\to \Gau(P)$ with contractible $\Gaut(P)$.
 We define an action of $\String$ on $\Gaut(P)$ such that $\Gaut(P)\to \String$
 is a smooth crossed module. Crossed modules are a source for Lie 2-groups
 (Example \ref{ex:bu1}) and in that way we obtain a Lie 2-group $\STRING$.
 \begin{nntheorem}[Theorem
  \ref{mainthm2}] $\STRING$ is a Lie 2-group model in the sense of Definition
  \ref{def:twomodel}.
 \end{nntheorem}
 The proof of this theorem relies on a comparison of the model $\String$ with
 the geometric realization of $\STRING$. Moreover, this direct comparison
 allows us to derive a comparison between the corresponding bundle theories and
 string structures, see Section \ref{sect:string_structures_and_comparison}.
 This explicit comparison is a distinct feature of our 2-group model that is
 not available for the other 2-group models. We show at the end of Section
 \ref{sect:string_structures_and_comparison} that any two string 2-group models
 (for instance the one from
 \cite{BaezCransStevensonSchreiber07From-loop-groups-to-2-groups} and from
 Section \ref{sect:the_string_group_as_a_2_group}) are comparable in the sense
 of the following result. This result is a variation of the uniqueness result
 from
 \cite{Schommer-Pries10Central-Extensions-of-Smooth-2-Groups-and-a-Finite-Dimensional-String-2-Group}
 for strict, infinite dimensional models of the string group, which is obtained
 by extracting a finite-dimensional presentation for the group stack associated
 to each Lie 2-group model in our setting.

 \begin{nntheorem}[Theorem
  \ref{thm:comparison}] If $\calg$ and $\calg'$ are smooth 2-group models for
  the string group, then there exists another smooth 2-group model $\calh$ and
  a span of morphisms
  \begin{equation*}
   \calg \xleftarrow{} \calh \xrightarrow{} \calg'
  \end{equation*}
  of smooth 2-group models.
 \end{nntheorem}
 
 Together with the other results of Section
 \ref{sect:string_structures_and_comparison} this then allows for an overall
 comparison between any kind of string bundle whose structure group is either a
 1-group or a 2-group.
 
 In an appendix we have collected some elementary facts about infinite
 dimensional manifolds and Lie groups. A second appendix gives a useful
 characterization of smooth weak equivalences between Lie 2-groups.
 
\end{tabsection}

\paragraph{Acknowledgements}

We thank Chris Schommer-Pries for pointing us to Stolz' string group model
\cite{Stolz}. We also thank Helge Gl{\"o}ckner, Eckhard
Meinrenken, Behrang Noohi, Friedrich Wagemann Danny Stevenson for discussions
and hints about various aspects. In particular, we want to thank Bas Janssens and Karl-Hermann Neeb for help with Remark \ref{rem:construction-of-basic-puH-bundle}. Finally, we thank Urs Schreiber and Christoph
Schweigert for comments on the draft.

While working on this paper TN was funded by the cluster of excellence
``Connecting particles with the cosmos'' and CW and CS were supported by the SFB 676
``Particles, Strings, and the Early Universe'' and CW was also supported by the GK 1670 ``Mathematics
Inspired by String Theory and QFT''. The paper was initiated during the
workshop ``Higher Gauge Theory, TQFT and Quantum Gravity'' in Lisbon 2011. All three authors thank the  Instituto Superior Técnico and the organizers for 
making it possible to attend this event.

\section{Preliminaries on gauge groups}
\label{sect:gaugegrp_lie}

\begin{tabsection}
 Throughout this paper Lie groups are permitted to be
 infinite-dimensional. More precisely, a Lie group is a group, together
 with the structure of a locally convex manifold such that the group
 operations are smooth, see Appendix \ref{app:locally_convex_manifolds}.
 The term topological group throughout refers to a group in the category of compactly
 generated weak Hausdorff spaces.

 Throughout this section we will use the following notation:
 \begin{itemize}
  \item $M$ is a compact manifold.
  \item $K$ is a metrizable Banach--Lie group (or equivalently a
        paracompact Banach--Lie group).
  \item $P$ is a smooth principal $K$-bundle over $M$.
 \end{itemize}
 
 Note that if $P$ is only a continuous principal bundle, then we can always
 find a smooth principal bundle which is equivalent to it \cite{wockel2}.
\end{tabsection}

\begin{definition}
 The group $\Aut(P)$ denotes the group of $K$-equivariant
 diffeomorphisms $f\from P\to P$. Identifying $M$ with $P/K$ we have a
 natural homomorphism
 \begin{equation*}
  Q\from \Aut(P)\to \Diff(M),\quad Q(f)([p])=[f(p)]
 \end{equation*}
 and we define the \emph{gauge group} by $\Gau(P):=\ker(Q)$.
\end{definition}

\begin{tabsection}
 It will be convenient to identify $\Gau(P)$ with $C^{\infty}(P,K)^K$,
 the smooth $K$-equivariant maps $P\to K$, via
 \[
  C^{\infty}(P,K)^K\ni f\mapsto (p\mapsto p\cdot f(p))\in\Gau(P).
 \]
 If $P$ is topologically trivial, then the left hand side
 $C^{\infty}(P,K)^K$ is isomorphic to $C^{\infty}(M,K)$. In
 \cite{wockel1} it is shown that in a certain sense this remains valid
 if $P$ is only locally trivial

 \begin{proposition}\label{prop:gauge}
  The group $\Gau(P) \cong C^{\infty}(P,K)^K$ admits the structure of a
  Fr{\'e}chet Lie-group modeled on the gauge algebra
  $\gau(P):= C^{\infty}(P,\fk)^K$ of smooth equivariant maps
  $P\to\fk$. If $\exp\from \fk\to K$ is the exponential function of $K$,
  then
  \begin{equation}\label{eqn:exponential_function_for_gauge_group}
   \exp_{*}\from C^{\infty}(P,\fk)^K \to C^{\infty}(P,K)^K,\quad \xi \mapsto
   \exp\circ\xi
  \end{equation}
  is an exponential function and a local diffeomorphism.
 \end{proposition}
 
\end{tabsection}

\begin{proof}
 The proof of this proposition can be found in \cite[Theorem 1.11 and
 Lemma 1.14(c)]{wockel1}. We will therefore only sketch the arguments
 that become important in the sequel. \\
 
 Let $N$ be a manifold with boundary (the boundary might be empty)
 modeled on a locally convex space. The space $C^{\infty}(N,K)$ can be
 given a topology by pulling back the compact open topology along
 \[
  C^{\infty}(N,K)\to\prod_{i=0}^{\infty} C^0(T^iN,T^iK)
 \]
 where $T^iN$ denotes the $i$-th iterated tangent bundle.
 We refer to this topology as the $C^{\infty}$-topology. This also
 applies to the Lie algebra $\fk$ of $K$ and induces a
 locally convex vector space topology on $C^{\infty}(N,\fk)$. Moreover,
 $C^{\infty}(N,\fk)$ is a Fr{\'e}chet space if $N$ is finite-dimensional
 \cite{Glockner02}. If we now restrict to the case where $N$ is compact and if
 $\varphi\from U\subset K\to W\subset \fk$ is a chart satisfying $\varphi(e)=0$,
 then $C^{\infty}(N,W)$ is in particular open in $C^{\infty}(N,\fk)$ and
 thus
 \begin{equation}\label{eqn:chart_for_mappting_group}
  \varphi_{*}\from C^{\infty}(N,U)\to C^{\infty}(N,W),\quad 
  \gamma\mapsto \varphi\circ \gamma	
 \end{equation}
 defines a manifold structure on $C^{\infty}(N,U)$. It can be shown that
 the (point-wise) group structures are compatible with this smooth
 structure and that it may be extended to a Lie group structure on
 $C^{\infty}(N,K)$. 
 Details of this construction can be found in
 \cite{wockel3} and
 \cite{GlocknerNeeb08Infinite-dimensional-Lie-groups-I}.

 The aforementioned topologies also endow the subspaces $C^{\infty}(P,K)^{K}$ and
 $C^{\infty}(P,\fk)^{K}$ with the structure of topological groups and
 $C^{\infty}(P,\fk)^{K}$ with the structure of a topological Lie
 algebra, both with respect to point-wise operations. The exponential
 function $\exp\from\fk\to K$ is $K$-equivariant
 and, by the inverse function theorem for Banach spaces, a local
 diffeomorphism. It thus defines in particular a map
 \begin{equation*}
  \exp_{*}\from C^{\infty}(P,\fk)^K \to C^{\infty}(P,K)^K,\quad \xi \mapsto
  \exp\circ\xi
 \end{equation*}
 Just as in the case of a compact manifold with boundary $N$, it can be
 shown that this map restricts to a bijection on some open subset of
 $C^{\infty}(P,\fk)^K$, which then gives rise to a manifold structure
 around the identity in $C^{\infty}(P,K)^K$ that can be enlarged to a
 Lie group structure. The details of this are spelled out in
 \cite[Propositions 1.4 and 1.8]{wockel1}.
\end{proof}

\begin{lemma}\label{lem:gau_is_paracompact}
 The topology underlying $\Gau(P)$ is metrizable.
\end{lemma}

\begin{proof}
 We first note that $C^{\infty}(N, K)$ is metrizable for
 finite-dimensional $N$ since $C^{0}(T^{i}N,T^{i} K)$ is so
 \cite[X.3.3]{Bourbaki98General-topology} and countable products of
 metrizable spaces are metrizable. From \cite[Proposition 1.8]{wockel1}
 it follows that $\Gau(P)$ is identified with a closed subspace of
 $C^{\infty}(\coprod \ol{V_{i}},K)$, where $V_{i},...,V_{n}$ is a cover
 of $M$ such that $\ol{V_{i}}$ is a manifold with boundary and
 $\left.P\right|_{\ol{V_{i}}}$ is trivial. Since
 $C^{\infty}(\coprod \ol{V_{i}},K)$ is metrizable, $\Gau(P)$ is so
 as well.
\end{proof}

\begin{remark}
 (\cite[Remark 1.18]{wockel1}) There also is a continuous version of the
 gauge group, namely the group of $K$-equivariant homeomorphisms $P\to P$
 covering the identity on $M$. This group will be denoted $\Gauc(P)$. As
 above, we have that $\Gauc(P)\cong C(P,K)^{K}$ and since $C(X,K)$ is a
 Lie group modeled on $C(X,\fk)$ for each compact topological space $X$
 (with respect to the compact-open topology, cf.\
 \cite{GlocknerNeeb08Infinite-dimensional-Lie-groups-I}) the above proof
 carries over to show that $\Gauc(P)$ is also a metrizable Lie group
 modeled on $C(P,\fk)^{K}$.
\end{remark}

Now \cite[Proposition 1.20]{wockel1} and Theorem \ref{thm:CW-type} imply

\begin{proposition}\label{prop:gau_homotopy_equivalence}
 The canonical inclusion
 \begin{equation}\label{eqn:hthy_equiv_gauge}
  \Gau(P)\hookrightarrow \Gauc(P).
 \end{equation}
 is a homotopy equivalence.
\end{proposition}

\begin{tabsection}
 In the sequel we will also need the following slight variation. Consider a central extension
 \begin{equation*}
  Z\to \widehat{K}\to K
 \end{equation*}
 of Banach--Lie groups admitting smooth local sections. Similar to
 $\Gau(P)\cong C^{\infty}(P,K)^{K}$, the groups $C^{\infty}(G,Z)$ and
 $C^{\infty}(P,\widehat{K})^{K}$ possess Lie group structures,
 modeled on $C^{\infty}(G,\mathfrak{z})$ and
 $C^{\infty}(P,\widehat{\mathfrak{k}})^{K}$ \cite[Appendix
 A]{NeebWockel07}, \cite[Theorem 1.11]{wockel1}. As in Proposition
 \ref{prop:gauge}, charts can be obtained from the exponential map
 \begin{equation*}
  \exp_{*}\from C^{\infty}(P,\widehat{\mathfrak{k}})^{K}\to C^{\infty}(P,\widehat{K})^{K},\quad \xi\mapsto \exp \op{\circ} \xi.
 \end{equation*}
 Moreover this is a central extension, as we show in proposition \ref{prop:local_sections}.
\end{tabsection}

\begin{lemma}\label{lem:metrizable_fiber_product}
 (\cite{EtterGriffin54}) If $F\to E\to B$ is a fiber bundle with $F$ and
 $B$ metrizable, then $E$ is metrizable.
\end{lemma}

\begin{proposition}\label{prop:local_sections}
 Let $Z\to \widehat{K}\xrightarrow{q} K$ be a central extension of
 Banach--Lie groups, admitting smooth local sections. Then the exact
 sequence of Fr\'echet--Lie groups
 \begin{equation}\label{eqn:central_extension_of_Gau(P)}
  C^{\infty}(M,Z)\to C^{\infty}(P,\widehat{K})^{K}\to C^{\infty}(P,K)^{K}
 \end{equation}
 also admits smooth local sections. Moreover, $C^{\infty}(M,\widehat{K})^{K}$ is
 metrizable if $Z$ and $K$ are so.
\end{proposition}

\begin{proof}
 We have to recall some facts on the construction of the Lie group
 structure from \cite[Appendix A]{NeebWockel07} and \cite[Proposition
 1.11]{wockel1}. Let $V_{1},...,V_{n}$ be an open cover of $G$ such that
 each $\ol{V_{i}}$ is a manifold (with boundary) and such that there exist smooth
 sections $\sigma_{i}\from \ol{V_{i}}\to P$. These give rise to smooth
 transition functions $k_{ij}\from\ol{V_{i}}\cap \ol{V_{j}}\to K$ and we
 have that
 \begin{equation*}
  \gamma \mapsto \Sigma(\gamma):= (\gamma \op{\circ} \sigma_{i} )_{i=1,...,n}
 \end{equation*}
 induces an isomorphism
 \begin{equation*}
  C^{\infty}(P,K)^{K}\cong         
  \{(\gamma_{i})_{i=1,...,n}\in\prod_{i=1}^{n} C^{\infty}(\ol{V_{i}},K)         
  \mid \gamma_{i}=k_{ij}\cdot \gamma_{j}\cdot k_{ji}         
  \text{ on }\ol{V_{i}}\cap\ol{V_{j}} \}         
 \end{equation*}
 If now $\exp\from\mathfrak{k}\to K$ restricts to a diffeomorphism
 $\exp\from W\to U$, then we have that
 \begin{equation*}
  \mathfrak{W}:=  \{(\gamma_{i})_{i=1,...,n}\in\prod_{i=1}^{n} C^{\infty}(\ol{V_{i}},W)         
  \mid \gamma_{i}=k_{ij}\cdot \gamma_{j}\cdot k_{ji}
  \text{ on }\ol{V_{i}}\cap\ol{V_{j}} \}         
 \end{equation*}
 maps under $\Sigma^{-1}$ to a neighborhood
 $\Sigma^{-1}(\mathfrak{W})$ of the identity on which $\exp_{*}$ restricts to a
 diffeomorphism (cf.\ \cite[Proposition 1.11]{wockel1}). Note that we
 may also assume without loss of generality that there exists a smooth section
 $\tau\from U\to \widehat{K}$ of $q$ satisfying
 $\tau(1_{K})=1_{\widehat{K}}$.

 Next we choose a smooth partition of unity
 $\lambda_{i}\from V_{i}\to [0,1]$. For
 $\gamma\in \Sigma^{-1}( \mathfrak{W})$ we then set
 \begin{equation*}
  \Lambda _{i}(\gamma):=\exp_{*}(\sum_{j\leq i}\lambda_{i}\cdot \log_{*}(\gamma))\cdot
  \exp_{*}(\sum_{j< i}\lambda_{i}\cdot \log_{*}(\gamma))^{-1}
 \end{equation*}
 and note that we have
 \begin{equation*}
  \gamma=\Lambda_{n}( \gamma)\cdot \Lambda_{n-1}( \gamma)\cdots \Lambda_{1}(\gamma).
 \end{equation*}
 Moreover, $\lambda_{i}(\pi(p))=0$ implies $\Lambda_{i}( \gamma)(p)=1$
 and thus $\op{supp}(\Lambda_{i}( \gamma))\subset V_{i}$. Moreover, we have
 $\Sigma(\Lambda_{i}(\gamma))_{i}\in C^{\infty}(\ol{V_{i}},W)$ by the
 definition of $\mathfrak{W}$.

 We now use all the data that we collected so far to define lifts of
 each $\Lambda_{i}(\gamma)$. To this end we first introduce functions
 $k_{i}\from \left.P\right|_{\overline{V}_{i}}\to K$, defined by
 $p=\sigma_{i}(\pi(p)).k_{i}(p)$. Then the assignment
 \begin{equation}\label{eqn:lift_construction_1}
  \left.P\right|_{V_{i}}\ni p\mapsto k_{i}(p). \tau\left(\Sigma(\Lambda_{i}(\gamma))_{i}(\pi(p))\right)
 \end{equation}
 is smooth since $\tau$ and $\Sigma(\Lambda_{i}(\gamma))_{i}$ are so and
 equivariant since $k_{i}$ is so. Moreover,
 \eqref{eqn:lift_construction_1} vanishes on a neighborhood of each
 point in $\partial \ol{V_{i}}$ since $\lambda_{i}$ and thus
 $\tau \op{\circ} \Sigma(\Lambda_{i}(\gamma))_{i}$ do so. Consequently,
 we may extend \eqref{eqn:lift_construction_1} by $e_{\widehat{K}}$ to
 all of $P$, defining  a lift $\Theta_{i}(\gamma)$ of
 $\Lambda_{i}(\gamma)$. Indeed, we have for $p\in \pi^{-1}(V_{i})$
 \begin{multline*}
  q(\Theta_{i}(\gamma)(p))=
  q\left(k_{i}(p).
  \tau\left(\Sigma(\Lambda_{i}(\gamma))_{i}(\pi(p))\right)\right)=
  k_{i}(p).q\left(\tau(\Sigma(\Lambda_{i}(\gamma))_{i}(\pi(p)))\right)=\\
  k_{i}(p).\Sigma(\Lambda_{i}(\gamma))_{i}(\pi(p))=k_{i}(p).\Lambda_{i}(\gamma)(\sigma_{i}(\pi(p))) =\Lambda_{i}(\sigma_{i}(\pi(p)).k_{i}(p))=\Lambda_{i}(\gamma)(p)
 \end{multline*}
 and for $p\notin \pi^{-1}(V_{i})$ we have
 $q(\Theta_{i}(\gamma)(p))=q(e_{\widehat{K}})=e_{K}=\Lambda_{i}(\gamma)(p)$.
 Eventually,
 \begin{equation*}
  \Theta(\gamma):=\Theta_{n}(\gamma)\cdot \Theta_{n-1}(\gamma)\cdots \Theta_{1}(\gamma)
 \end{equation*}
 defines a lift of $\gamma$, since we have
 \begin{multline*}
  q_{*}(\Theta_{n}(\gamma)\cdot \Theta_{n-1}(\gamma)\cdots \Theta_{1}(\gamma))=
  q_{*}(\Theta_{n}(\gamma))\cdot q_{*}(\Theta_{n-1}(\gamma))\cdots q_{*}(\Theta_{1}(\gamma))=\\
  \Lambda_{n}(\gamma)\cdots \Lambda_{n-1}(\gamma)\cdots \Lambda_{1}(\gamma)=\gamma.
 \end{multline*}
 Since $\Theta_{i}(\gamma)$ is constructed in terms of push-forwards of
 smooth maps, it depends smoothly on $\gamma$ and so does
 $\Theta(\gamma)$.

 The previous argument shows in particular that
 \eqref{eqn:central_extension_of_Gau(P)} is a fiber bundle (cf.\
 \ref{lem:submersion}). As in Lemma \ref{lem:gau_is_paracompact} one
 sees that $C^{\infty}(M,Z)$ is metrizable if $Z$ is so, and thus the
 last claim follows from Lemma \ref{lem:metrizable_fiber_product}.
\end{proof}

\begin{remark}\label{rem:generality1}
 Note that all results of this section remain valid in more general
 situations. For instance, if we replace $K$ by an arbitrary Lie group
 with exponential function that is a local diffeomorphism, then
 $\Gau(P)$ is a Lie group, modeled on $\gau(P)$. Moreover,
 \eqref{eqn:exponential_function_for_gauge_group} still defines an
 exponential function which itself is a local diffeomorphism. If, in
 addition, $K$ is metrizable, then the proof of Lemma
 \ref{lem:gau_is_paracompact} shows that $\Gau(P)$ is also
 metrizable.
 
 Proposition \ref{prop:local_sections} generalizes to the situation
 where $Z\to \widehat{K}\to K$ is a central extension of Lie groups for
 which $\widehat {K}$ and $K$ have exponential functions that are local
 diffeomorphisms. Since its proof only uses the fact that
 $\widehat{K}\to K$ has smooth local sections,
 \eqref{eqn:central_extension_of_Gau(P)} still admits smooth local
 sections in this case.

 This shows in particular that the construction applies to the smooth
 principal bundle $\Omega G\to P_{e}G \to G$, where $\Omega G$ denotes the
 group of based smooth loops (as for instance in
 \cite[Section 3]{BaezCransStevensonSchreiber07From-loop-groups-to-2-groups}) and
 the universal central extension $U(1)\to \widehat{\Omega G}\to \Omega G$.
\end{remark}

\section{The string group as a smooth extension of \texorpdfstring{$G$}{G}}\label{sec:group}

In this section we want to give a smooth model for the string group. Note that Lie groups in our
setting are modeled on arbitrary locally convex spaces (cf.\ Appendix \ref{app:locally_convex_manifolds}),
so they may be in particular infinite-dimensional. Our construction is based on 
\cite[Section 5]{Stolz}. We are mainly interested in the case $G = Spin(n)$ but we define more generally:

\begin{definition}
 \label{def:model} Let $G$ be a compact, simple and 1-connected Lie group. A
 \emph{smooth string group model} for $G$ is a Lie group $\widehat{G}$ together
 with a smooth homomorphism
 \begin{equation*}
  \widehat{G}\xrightarrow{q} G
 \end{equation*}
 such that $q$ is a Serre fibration,
 $\pi_{k}(\widehat{G})=0$ for $k\leq 3$ and that $\pi_{i}(q)$ is an isomorphism
 for $i> 3$.
\end{definition}

\begin{proposition}[Cartan \cite{Car36}]
Let $G$ be a compact, simple and 1-connected Lie group. Then
$$ \pi_2(G) = 0 \quad \text{and} \quad \pi_3(G) \cong H^3(G,\Z) \cong \Z\text{.}$$
\end{proposition}

\begin{corollary}\label{cor:string-is-not-finite-dimensional}
 If $\widehat{G}\xrightarrow{q} G$ is a smooth string group model, then
 \begin{enumerate}
        \renewcommand{\labelenumi}{\theenumi}
        \renewcommand{\theenumi}{\arabic{enumi}.}
  \item \label{item:K(Z,2)_1} $\ker(q)$ is a $K(\Z,2)$ (i.e.,
        $\pi_{k}(\ker(q))\cong \Z$ for $k=2$ and vanishes for $k\neq 2$);
  \item $\widehat{G}$ cannot be finite-dimensional.
 \end{enumerate}
\end{corollary}

\begin{proof}
 \begin{enumerate}
  \item This follows from the long exact homotopy sequence.
  \item If $\widehat{G}$ were finite-dimensional, then it would have $\ker(q)$
        as a closed Lie subgroup. But by \ref{item:K(Z,2)_1} we have
        ${H}^{2n}(\ker(q),\Z)\cong H^{2n}(K(\Z,2),\Z)\cong\Z$, a contradiction
        (one could also use \cite[Theorem
        2]{GotayLashofSniatyckiWeinstein83Closed-forms-on-symplectic-fibre-bundles}
        to conclude that $H^{2}(\widehat{G})=0$ and $H^{2}(\ker(q))\neq 0$ implies
        that $\widehat{G}$ has infinite cohomology).
 \end{enumerate}
\end{proof}

\begin{remark}
 In the definition of smooth string group model (Definition \ref{def:model}) it
 is possible to impose other conditions on $q$.
 \begin{enumerate}
  \item A seemingly stronger one is to require $q$ to be a topological locally
        trivial bundle. However, in our particular situation $q$ is a Serre
        fibration if and only if it is a topological locally trivial bundle. In
        fact, let $U\subset G$ be a contractible neighborhood of $g\in G$. Then
        the cofibration $\{g\}\hookrightarrow U$ allows a lift
        $\sigma\from U\to \widehat{G}$ since trivial cofibrations have the left
        lifting property for Serre fibrations. This then implies that $q$ is a
        locally trivial bundle, since a continuous section
        $\sigma\from U\to \widehat{G}$ yields the trivialization
        \begin{equation*}
         q^{-1}(U)\to U\times\ker(q),\quad g\mapsto (q(g),\sigma(q(g))^{-1}\cdot g).
        \end{equation*}
        On the other hand, if $q$ is a locally trivial bundle, then it is a
        Hurewicz fibration by \cite[Corollary
        2.8.14]{Spanier66Algebraic-topology} and thus in particular a Serre
        fibration.
  \item Another possibility is to drop the condition that the map
        $q: \widehat{G} \to G$ is a Serre-fibration. From the point of homotopy
        theory this might even be more natural. However in this case the kernel
        of $q$ is not a $K(\Z,2)$ anymore but only the homotopy kernel. The
        second assertion of Corollary
        \ref{cor:string-is-not-finite-dimensional} would still remain valid in
        this case. To show this, it suffices to show that the cohomology of
        $\widehat G$ is non-zero in infinitely many degrees. This can be done
        using the Serre spectral sequence for the fibre sequence
        $K(\Z,2) \to \widehat{G} \to G$. More precisely if
        $x \in H^2(K(\Z,2),\Z)$ and $y \in H^3(G,\Z)$ are generators of the
        respective cohomology groups, one can show that the elements
        $x^{n} y \in H^3\big(G,H^{2n}(K(\Z,2),\Z)\big)$ for even $n$ survive
        until the $E_\infty$-term and thus give rise to non-vanishing classes
        in $H^{2n+3}(\widehat G,\Z)$.
        
 \end{enumerate}
 
\end{remark}

Now we come to the construction of our string group model. Let $\calh$ be an infinite-dimensional separable Hilbert space. Then it is well known that the projective unitary group $PU(\calh)$, together with the norm topology is a $K(\Z,2)$  \cite{Kuiper65The-homotopy-type-of-the-unitary-group-of-Hilbert-space}, so
that $BPU(\calh)$ is a $K(\Z,3)$. This induces a bijection between isomorphism classes of $PU(\calh)$-bundles over a manifold $M$ and $H^3(M,\Z)$. 

Now there is a canonical generator $1 \in H^3(G,\Z)$. Let $P\to G$ be a principal
$PU(\calh)$-bundle over $G$ that represents this generator. Note that $PU(\calh)$
is a Banach--Lie group (see
\cite{GlocknerNeeb03Banach-Lie-quotients-enlargibility-and-universal-complexifications}
and references therein) which is paracompact by \cite[Theorem VIII.2.4]{Dugundji66Topology} and \cite[Theorem I.3.1]{Bredon72Introduction-to-compact-transformation-groups}. In particular, it is metrizable.
We can choose $P$ to be smooth \cite{wockel2} and apply the results from Section \ref{sect:gaugegrp_lie}  (we will discuss in Remark \ref{rem:construction-of-basic-puH-bundle} the problem of giving a geometric construction of this bundle). Recall in particular the map
$$ Q:\Aut(P) \to  \Diff(G) $$
that sends a bundle automorphism to its underlying diffeomorphism of the base.

\begin{definition}\label{def:string}
Let $G$ be a compact, simple and 1-connected Lie group and $P \to G$ be a principal $PU(\cal{H})$-bundle representing the generator $1 \in H^3(G,\Z)$. Then we set
$$ \String := \{ f\in\Aut(P) \mid Q(f)\in G \subset \Diff(G) \} $$
where the inclusion $G \hookrightarrow \Diff(G)$ sends $g$ to left multiplication with $g$. In other words: $\String$ is the group consisting of bundle automorphisms that cover left multiplication in $G$.
\end{definition}

\begin{tabsection}
 Note that there is also a continuous version of $\String$, given by
 \begin{equation*}
  \Stringc:=\{f\in\Homeo(P)\mid f\text{ is } PU(\calh) \text{-equivariant and }Q(f)\in G\subset\Diff(G) \}.
 \end{equation*}
 The motivation for constructing a smooth model for the String group as
 in the present paper now comes from the following fact \cite{Stolz}.
 For the sake of completeness we include (a part of) the proof here.
\end{tabsection}

\begin{proposition}[Stolz]\label{prop:stolz}
The fibration $Q: \Stringc \to G$ is a 3-connected cover of $G$, i.e. $\pi_i(\Stringc) = 0$ for $i \leq 3$ and $\pi_i(Q)$ is an isomorphism for $i > 3$.
\end{proposition}
\begin{proof}
Pick a point $p \in P$ in the fiber over $1 \in G$. Let $ev$ be the evaluation that sends a bundle automorphism $f$ to $f(p)$.
Then we obtain a diagram
$$
\xymatrix{
{\Gauc(P)} \ar[r] \ar[d]^{ev} & {\Stringc} \ar[d]^{ev} \ar[r]^Q & G \ar[d]_{\mathrm{id}}\\
PU(\calh) \ar[r] & P  \ar[r]_\pi & G 
}
$$
Now  \cite[Lemma 5.6]{Stolz} asserts that $ev: \Gauc(P) \to PU(\calh)$ 
is a (weak) homotopy equivalence. The long exact homotopy sequence and the Five Lemma then 
show that then $ev: \Stringc \to P$ is also a homotopy equivalence. Hence it remains to show that $P \to G$ is a 3-connected cover. By definition of $P$ its classifying map 
$$ p: G \longrightarrow BPU(\calh) \simeq K(\Z,3) $$
is a generator of $H^3(G,\Z)$, hence it induces isomorphisms on the first three homotopy groups. Thus the pullback $P \cong p^*EPU(\calh)$ of the contractible space $EPU(\calh)$ kills exactly the first three homotopy groups, i.e. $P$ is a 3-connected cover.
\end{proof}

In the rest of this section we want to prove the following modification and
enhancement of the preceding proposition. For its formulation recall that an
extension of Lie groups is a sequence of Lie groups $A\to B\to C$ such that
$B$ is a smooth locally trivial principal $A$-bundle over $C$
\cite{Neeb07}.

\begin{theorem}\label{mainthm1}
 $\String$ is a smooth string group model according to Definition
 \ref{def:model}. Moreover, $\String$ is metrizable and there exists a
 Fr\'echet--Lie group structure on $\String$, unique up to isomorphism,
 such that
 \begin{equation}\label{eqn:sting_extension}
  \Gau(P)\to \String \to G
 \end{equation}
 is an extension of Lie groups.
\end{theorem}

\begin{proof}
 We first show existence of the Lie group structure. To this end we
 recall that there exists an extension of Fr\'echet--Lie groups
 \begin{equation}\label{eqn:Aut_extension}
  \Gau(P)\to \Aut(P)_{0}\to \Diff(G)_{0},
 \end{equation}
 where $\Aut(P)_{0}$ is the inverse image $Q^{-1}(\Diff(G)_{0})$ of the
 the identity component $\Diff(M)_{0}$ \cite[Theorem 2.14]{wockel1}. The
 embedding $G\hookrightarrow \Diff(G)_{0}$ given by left translation
 gives by the exponential law
 \cite{GlocknerNeeb08Infinite-dimensional-Lie-groups-I} a smooth
 homomorphism of Lie groups since the multiplication map
 $G\times G\to G$ is smooth. Pulling back \eqref{eqn:Aut_extension}
 along this embedding then yields the extension
 \eqref{eqn:sting_extension}. Moreover, $\String$ is metrizable by Lemma
 \ref{lem:gau_is_paracompact} and Lemma
 \ref{lem:metrizable_fiber_product}.

 We now discuss the uniqueness assertion, so let
 $\Gau(P)\to H_{i} \xrightarrow{q_i} G$ for $i=1,2$ be two extensions of
 Lie groups. The requirement for it to be a locally trivial smooth
 principal bundle is equivalent to the existence of a smooth local
 section of $q_i$ and we thus obtain a derived extension of Lie algebras
 \begin{equation*}
  \gau(P)\to L(H_{i}) \xrightarrow{L(q_i)} \fg.
 \end{equation*}
 The differential of the local smooth section implements a linear
 continuous section of $L(q_i)$ and thus we have a (non-abelian)
 extension of Lie algebras in the sense of \cite{Neeb06}. Now the
 equivalence classes of such extensions are parametrized by
 $H^{2}(\fg,\mathfrak{z}(\gau(P)))$ \cite[Theorem II.7]{Neeb06}. Since
 $\gau(P)= C^{\infty}(P,\mathfrak{pu}(H))^{K}$ we clearly have
 $\mathfrak{z}(\gau(P))=C^{\infty}(P,\mathfrak{z}(\mathfrak{pu}(H)))^{K}$,
 which is trivial since $\mathfrak{z}(\mathfrak{pu}(H))$ is so.
 Consequently, we have a morphism
 \begin{equation*}
  \xymatrix{
  \gau(P)\ar[r]\ar@{=}[d]&L(H_{1})\ar[r]\ar[d]^{\varphi}&\fg\ar@{=}[d]\\
  \gau(P)\ar[r]&L(H_{2})\ar[r]&\fg
  }
 \end{equation*}
 of extensions of Lie algebras. The long exact homotopy sequence for the
 fibration $\Gau(P)\to H_{i} \xrightarrow{q} G$ shows that $H_{i}$ is
 1-connected, and so $\varphi$ integrates to a morphism
 \begin{equation*}
  \vcenter{  \xymatrix{
  \Gau(P)\ar[r]\ar@{=}[d]&H_{1}\ar[r]\ar[d]^{\Phi}&G\ar@{=}[d]\\
  \Gau(P)\ar[r]&H_{2}\ar[r]& G 
  }}
 \end{equation*}
 of Lie groups (note that $H_{i}$ is regular by \cite[Theorem
 5.4]{OmoriMaedaYoshiokaKobayashi83On-regular-Frechet-Lie-groups.-V.-Several-basic-properties}
 since $\Gau(P)$ and $G$ are so). Since $\Phi$ makes this diagram
 commute it is automatically an isomorphism.

 It remains to show that $\String$ is a smooth model for the String
 group. We have the following commuting diagram $$\vcenter{\xymatrix{
 \Gau(P) \ar[r]\ar[d] & \String\ar[d] \ar[r] & G \ar@{=}[d]\\
 \Gauc(P) \ar[r] & \Stringc \ar[r] & G }}. $$ By Proposition
 \ref{prop:gau_homotopy_equivalence} the inclusion
 $\Gau(P) \hookrightarrow \Gauc(P)$ is a homotopy equivalence. Since,
 furthermore, $\String \to G$ and $\Stringc \to G$ are bundles, they are
 in particular fibrations and we obtain long exact sequences of homotopy
 groups.
 Applying the Five Lemma we see that the maps
 $\pi_n(\String) \to \pi_n(\Stringc)$ are isomorphisms for all $n$. By
 Proposition \ref{prop:stolz} we know that $\Stringc$ is a 3-connected
 cover, hence also $\String$.

\end{proof}

\begin{remark}
 Note that the proof of the uniqueness assertion only used the fact that
 the center of $\gau(P)$ is trivial. In fact, this shows that for an
 arbitrary (regular) Lie group $H$ which is a $K(\Z,2)$ and has trivial
 $\mathfrak{z}(L(H))$ there exists, up to isomorphism, at most one Lie group
 $\widehat{H}$, together with smooth maps $H\to \widehat{H}$ and
 $\widehat{H}\to G$ turning
 \begin{equation*}
  H\to\widehat{H}\to G
 \end{equation*}
 into an extension of Lie groups. Moreover, the proof shows that the
 uniqueness is not only up to isomorphism of Lie groups, but even up to
 isomorphism of \emph{extensions}.
\end{remark}

\begin{remark}\label{rem:construction-of-basic-puH-bundle}
 Since the accessibility of the smooth model $\String$ for the string group
 depends on the accessibility of the basic $PU(\calh)$-bundle $P\to G$ we
 shortly comment on how this might (or might not) be constructed. A commonly
 used construction is to take the smooth path-loop fibration
 \begin{equation*}
  P_{e}G:=\{\gamma\in C^{\infty}(\R,G):\gamma(x+1)\cdot \gamma(x)^{-1}\text{ is constant for all }x\in\R\text{ and } \gamma(0)=e\}
 \end{equation*}
 and then consider $P_{e}G\to G$, $\gamma\mapsto \gamma(1)$ as a principal
 $\Omega G$-bundle over $G$
 If one now takes a positive energy representation $\rho\from LG\to PU(\calh)$
 of level $1$ and restricts it to $\Omega G$, then the resulting associated
 $PU(\calh)$-bundle over $G$ also has level $1$, i.e., represents a generator
 in $[G,BPU(\calh)]\cong\Z$ (a similar construction also works for the full
 loop group $LG$ if one considers the principal $LG$-bundle over $G$ consisting
 of the quasi-periodic paths).
 
 This construction yields a \emph{continuous} principal bundle for the
 \emph{strong operator topology} on $PU(\calh)$ (in which it also is a
 $K(\Z,2)$), since then $\rho$ is a homomorphism of topological groups.
 However, $\rho$ cannot be smooth in the norm topology\footnote{In the strong
 topology $PU(\calh)$ is not a Lie group.} and thus one cannot obtain a smooth
 basic $PU(\calh)$-bundle this way. In fact, the projective representation of
 $\Omega G$ corresponds to a linear representation of the universal central
 extension $\widehat{\Omega G}$. The associated Lie algebra
 $\widehat{\Omega \fg}=\R\oplus_{\omega}\Omega \fg$ has a three dimensional
 Heisenberg algebra as subalgebra, for instance generated by an appropriately
 normalized element $X\in \fg$ and the elements
 \begin{equation*}
  X_{1}(t):=\sin(2 \pi t)\cdot X\quad\text{ and }\quad X_{2}(t):=-\cos(2 \pi t)\cdot X+X
 \end{equation*}
 of ${\Omega \fg}\subset \widehat{\Omega \fg}$. Thus any strongly continuous
 representation of $\widehat{\Omega G}$ induces a strongly continuous
 2-parameter family $(W_{s,t})_{s,t\in \R}$ of operators satisfying the Weyl
 relations. Now von Neumann's Uniqueness Theorem implies that the
 representation is the tensor product of the standard representation of the CCR
 algebra and a trivial representation. If now the center of
 $\widehat{\Omega G}$ does not act by $0$, then this standard representation is
 not bounded. However, boundedness of a representation is a necessary condition
 for continuity (in the norm topology) and thus for smoothness. Consequently,
 any representation of $\widehat{\Omega G}$ or $\widehat{LG}$ for which the
 center does not act by $0$ cannot be smooth.
\end{remark}

\begin{problem}
 The previous remark puts a severe constraint on the way how the smooth basic
 $PU(\calh)$-bundle $P\to G$ may be constructed. The abstract existence result
 from \cite{wockel2} will also not lead to an explicit construction. We thus
 consider it as an interesting problem to give such a
 construction in explicit, at best in geometric terms.
\end{problem}

\section{2-groups and 2-group models}
\label{sect:2groups}

One of the main problems with string group models is that they are not very tightly determined. In fact, the underlying space is just determined up to weak homotopy equivalence. This implies that the group structure can only be determined up to $A_\infty$-equivalence and the smooth structure is not determined at all. Part of the problem is that one does not in general have good control over the fiber of $\String \to G$, only the underlying homotopy type is determined to be a $K(\Z,2)$. 

This problem can be cured by using 2-group models. This setting allows to fix the fiber more tightly. In particular there is a nice model of $K(\Z,2)$ as a 2-group, see Example \ref{ex:bu1} below and weak equivalences of 2-groups are more restrictive than homotopy equivalences of their geometric realizations. Moreover, this setting
implies a strong uniqueness condition, as we point out in the end of the section.

We first want to recall quickly the definition and some elementary properties of 2-groups. We restrict our attention  to strict Lie 2-groups in this paper which for simplicity we just call Lie 2-groups.

\begin{definition}
 A (strict) Lie 2-group is a category $\calg$ such that the set of
 objects $\calg_0$ and the set of morphisms $\calg_1$ are Lie groups,
 all structure maps
 \begin{equation*}
  s,t: \calg_1 \to \calg_0 \quad i: \calg_0 \to \calg_1 \quad \text{and} \quad \circ: \calg_1 \times_{\calg_0} \calg_1 \to \calg_1 
 \end{equation*}
 are Lie group homomorphisms and $s,t$ are
 submersions\footnote{Submersion in the sense that it is locally a
 projection, see Appendix \ref{app:locally_convex_manifolds}}. In the
 case that $\calg_0$ and $\calg_1$ are metrizable, we call $\calg$ a
 metrizable Lie 2-group. A morphism between 2-groups is a functor
 $f: \calg \to \calg'$ that is a Lie group homomorphism on the level of
 objects and on the level of morphisms.
\end{definition}

One reason to consider 2-groups here is that they can serve as models for topological spaces by virtue of the following construction.

\begin{definition}
Let $\calg$ be a Lie 2-group. Then the nerve $N\calg$ of the category $\calg$ is a simplicial manifold by Proposition \ref{prop:existence_of_fiber_products}. Using this we define the \emph{geometric realization} of $\calg$ to be the geometric realization of the simplicial space $N\calg$, i.e., the coend
\begin{equation*}
\int^{[n] \in \Delta} (N\calg)_n \times \Delta^n = \bigsqcup_{n}~ (N\calg)_n \times \Delta^n ~/~\sim .
\end{equation*}
Note that the coend is taken in the category of compactly generated weak Hausdorff spaces. The 
\emph{fat geometric realization} is obtained in an exactly analogous way after restricting the
simplex category $\Delta$ to the subcategory $\Delta_+$ which has only strictly increasing maps
(called ``fat'' because here the relations given by degeneracy maps are not quotiented out).
\end{definition}

\begin{example}
 \phantomsection\label{ex:bu1}
 \begin{enumerate}
  \item Consider the category $\calb U(1)$ with one object and
        automorphisms given by the group $U(1)$. This is clearly a Lie
        2-group. The geometric realization $|\calb U(1)|$ is the
        classifying space $BU(1)$, hence a $K(\Z,2)$. The 2-group
        $\calb A$ exists moreover for each abelian Lie group $A$.
  \item If $G$ is an arbitrary Lie group, then it gives rise to a
        2-group by considering it as category with only identity
        morphisms. More precisely, in this case $\calg_{0}=\calg_{1}=G$
        and all structure maps are the identity.
  \item Let $K \xrightarrow{\partial} L$ be a smooth crossed module of
        groups (\cite[Definition 3.1]{Neeb07}). Then we can form a Lie
        2-group $\calg$ using the Lie groups $\calg_0 := L$ and
        $\calg_1 := K \rtimes L$ together with the smooth maps
        $s(k,l) = l$, $t(k,l) = \partial(k)l$, $i(l) = (1,l)$ and
        $(k,l) \circ (k',l') = (kk',l)$.

        Conversely, if we start with a Lie 2-group $\calg'$, then
        $\ker(s')$ is a Lie subgroup of $\calg'_{1}$ by Lemma
        \ref{lem:submersion} and we obtain a crossed module by setting
        $K:=\ker(s)$, $L:=\calg'_{0}$,
        $\partial:=\left.t\right|_{\ker(s)}$ and defining the (left)
        action to be $k.l:=i(l^{-1})\cdot k\cdot i(l)$.
 \end{enumerate}
\end{example}

\begin{lemma}\label{lem:good}
If $\calg$ is a metrizable Lie 2-group, then
\begin{enumerate}\setlength{\itemsep}{-0.5ex}
\item all spaces $N\calg_n$ have the homotopy type of a CW complex;
\item the nerve $N\calg$ is good, i.e. all degeneracies are closed cofibrations;
\item the nerve $N\calg$ is proper, i.e Reedy cofibrant as a simplicial space (with respect to the Str{\o}m model structure \cite{Strom});
\item the canonical map from the fat geometric realization $\Vert N\calg \Vert$ to the ordinary geometric realization $|\calg|$ is a homotopy equivalence;
\item the geometric realization $|\calg|$ has the homotopy type of a CW-complex.
\end{enumerate}
\end{lemma}
\begin{proof}
1) First note that all the spaces $(N\calg)_n$ are subspaces of $(\calg_{1})^{n}$ and thus are metrizable. Hence by Theorem \ref{thm:CW-type} they have the homotopy type of a CW-complex.

2) Again using the fact that all $(N\calg)_n$  are metrizable and \cite[Theorem 7]{pal66} we see that they are well-pointed in the sense that the basepoint inclusion is a closed cofibration. A statement of Roberts and Stevenson \cite[Proposition 23]{RSDraft} then shows that $N\calg$ is good, i.e., degeneracy maps are closed cofibrations. We roughly sketch a variant of their argument here: 
By the fact that $\calg$ is a 2-group we can write the nerve as 
\begin{equation*}
\xymatrix{
{\cdots}  
\ar@<1.3ex>[r]
\ar@<0.3ex>[r]
\ar@<-0.7ex>[r]
\ar@<-1.7ex>[r]
&
{\ker(s) \times \ker(s) \times \calg_0}
\ar@<0.9ex>[r]
\ar@<-0.1ex>[r]
\ar@<-1.1ex>[r]
&
{\ker(s) \times \calg_0 }
\ar@<0.3ex>[r]
\ar@<-0.7ex>[r]
&
\calg_0
}
\end{equation*}
where the decomposition is a decomposition on the level of topological spaces. Hence to show that the degeneracies are 
closed cofibrations it suffices to show that $\ker(s)$ is well-pointed. But it is a retract of $\calg_1 = \calg_0 \times \ker{s}$ hence well pointed by the fact that $\calg_1$ is well pointed. 

3) Now we know that $N\calg$ is good and in this case \cite[Corollary 2.4(b)]{Lewis} implies that $N\calg$ is also proper.

4) By \cite[Proposition A1]{Seg74} (resp \cite[Proposition 1]{Dieck74}) the fat and the ordinary geometric realizations are homotopy equivalent. 

5) Since all the spaces $(N\calg)_n$ have the homotopy type of a CW-complex, also the fat geometric realization has the homotopy type of a CW complex \cite[Proposition A1]{Seg74}. Thus also the ordinary realization by 4).
\end{proof}

\begin{proposition}\label{prop:equivreal}
If $\calg$ and $\calg'$ are metrizable Lie 2-groups and $f: \calg \to \calg'$ is a homomorphism that is a weak homotopy equivalence on objects and morphisms, then 
$$|f|: \quad |\calg| \to |\calg'|$$
is a homotopy equivalence.
\end{proposition}
\begin{proof}
First note that $Nf: N\calg \to N\calg'$ is a level-wise weak homotopy equivalence. For 
0-simplices and 1-simplices this is true by assumption and for the higher simplices it follows again from the product structure of the nerves given in the proof of Lemma \ref{lem:good} and the fact that $Nf$ is also a product map (note that $f$ induces a weak homotopy equivalence $\ker(s)\to\ker(s')$). Then using \cite[Proposition A4]{May74} and the fact that $N\calg$ and $N\calg'$ are proper we conclude that also $|f|: |\calg| \to |\calg'|$ is a weak homotopy equivalence. But since the geometric realizations have the homotopy type of a CW-complex, Whitehead's theorem shows that $|f|$ is an honest homotopy equivalence.
\end{proof}

For smooth groupoids there is a notion of weak equivalence which is equivalent to equivalence of the associated stacks, see e.g. \cite[Definition 58 and Proposition 60]{metzler}. We adopt this for 2-groups.

\begin{definition}
\label{def:swe}
A morphism $f: \calg \to \calg'$ of Lie 2-groups is called \emph{smooth weak equivalence} if the following conditions are satisfied:
\begin{enumerate}
\item 
it is smoothly essentially surjective: the map
\begin{equation*}
s \circ \mathrm{pr}_2: {\calg}_0 {~}_{f_0}\!\times_{t} {\calg'}_1 \to {\calg'}_0
\end{equation*}
is a surjective submersion.
\item
it is smoothly fully faithful: the diagram
\begin{equation*}
\xymatrix{
{\calg}_1 \ar[r]^{f_1} \ar[d]_{s \times t} & {\calg'}_1 \ar[d]^{s \times t} \\
{\calg}_0 \times {\calg}_0 \ar[r]_-{f_0 \times f_0} & {\calg'}_0 \times {\calg'}_0
}
\end{equation*}
is a pullback diagram.
\end{enumerate}

\end{definition}

\begin{proposition}\label{morita}
Let $f: \calg \to \calg'$ be a smooth weak equivalence between metrizable 2-groups. Then $|f|: |\calg| \to |\calg'|$ is a homotopy equivalence.
\end{proposition}
\begin{proof}
A smooth weak equivalence between 2-groups is in particular a topological weak equivalence of the underlying topological groupoids. But then a result of Gepner and Henriques \cite{GepHen} or Noohi
\cite[Theorem 6.3 and Theorem 8.2.]{noohi} implies that the induced morphism $\Vert f \Vert: \Vert \calg \Vert \to \Vert \calg' \Vert$ between the fat geometric realizations is a weak equivalence. Again by the fact that the fat realizations are homotopy equivalent to the geometric realizations this completes the proof.
\end{proof}

\begin{definition}
 If $\calg$ is a Lie 2-group, then we denote by $\upi_{0}\calg$ the
 group of  isomorphism classes of objects in $\calg$ and by
 $\upi_{1}\calg$ the group of automorphisms of $1\in \calg_{0}$. Note
 that $\upi_{1}\calg$ is abelian. We call $\calg$ \emph{smoothly
 separable} if $\upi_{1}\calg$ is a {split Lie subgroup}\footnote{Split Lie subgroup in the sense of Definition \ref{defslipt}}
 of $\calg_{1}$ and $\upi_{0}\calg$
 carries a Lie group structure such that
 $\calg_{0}\to \upi_{0}\calg$ is a submersion.
\end{definition}

\begin{proposition}\phantomsection\label{prop:separable}
\begin{enumerate}
\item
A morphism between smoothly separable Lie 2-groups is a smooth weak equivalence if and only if it induces Lie group isomorphisms on $\upi_0$ and $\upi_1$.
\item
For a metrizable, smoothly separable Lie 2-group $\calg$ the sequence
\begin{equation*}
|\calb\upi_1\calg| \to |\calg| \to \upi_0\calg
\end{equation*}
is a fiber sequence of topological groups. Moreover, the right hand map is a fiber bundle and the left map is a homotopy equivalence to its fiber.
\end{enumerate}
\end{proposition}
\begin{proof}
The first claim will be proved in Appendix \ref{app:lie2groups}. We thus show the second. Let us first consider the morphism $q: \calg \to \upi_0\calg$ of 2-groups 
where $\upi_0\calg$ is considered as a 2-group 
with only identity morphisms. Let $\mathcal{K}$ be the level-wise kernel of this map, i.e., $\calk_{0}=\ker(q_{0})$ and $\calk_{1}=\ker(q_{1})$. Since $q_{1}=q_{0}\op{\circ} s$ it is a submersion, $\calk_{0}$ and $\calk_{1}$ are Lie subgroups and $\calk$ is a metrizable Lie 2-group.  Then $N\mathcal{K} \to N\calg \to N\upi_0\calg$ is an exact 
sequence of simplicial groups. It is easy to see that the geometric realization of this sequence is also exact, e.g., by using 
the fact that geometric realization preserves pullbacks \cite[Corollary 11.6]{May74}. Hence we have an exact sequence of topological groups.
\begin{equation*}
|\mathcal{K}| \to |\calg| \to \upi_0\calg
\end{equation*}
Moreover the right hand map is a $|\mathcal{K}|$-bundle since by the definition of smooth separability 
it admits local sections. Thus it only remains to show that $|\calb\upi_1\calg| \simeq |\mathcal{K}|$. 
Now the inclusion $\calb \upi_1\calg \to \mathcal{K}$ is a smooth weak equivalence, which we can either see using the first part of the Proposition or by a direct argument. Then Proposition \ref{morita} shows that the realization is a homotopy equivalence.
\end{proof}

\begin{definition}\label{def:twomodel}
 Let $G$ be a compact, simple and 1-connected Lie group. A
 \emph{smooth 2-group model} for the string group is a metrizable smooth
 2-group $\calg$ which is smoothly separable together with isomorphisms
 \begin{equation}\label{eqn:morphism}
  \upi_0\calg \iso G \qquad \text{and} \qquad \upi_1\calg \iso U(1)
 \end{equation}
 such that $ |\calg| \to G$ is a 3-connected cover.
\end{definition}

The following lemma is immediate and explains the last condition in the previous definition.

\begin{lemma}\label{lem:lifting}
 If $\calg$ is a smoothly separable 2-group such that $\upi_0\calg\cong G$
 and $\upi_1\calg \cong U(1)$, then the following are equivalent.
 \begin{enumerate}
     \renewcommand{\labelenumi}{\theenumi}
     \renewcommand{\theenumi}{\arabic{enumi}.}
  \item\label{item:lifting_1}  $\calg$ is a smooth 2-group model for the string group.
  \item\label{item:lifting_2} The surjective submersion $\calg_{0}\to G$ together with the $U(1)$-bundle 
        \begin{equation}\label{eqn:lifting_1}
        \calg_{1}\xrightarrow{s\times t} \calg_{0}\times_{G}\calg_{0}
        \end{equation}
        and the composition in the groupoid $\calg$ is a model for a basic bundle gerbe over $G$, i.e. a generator of $H^3(G,\Z) = \Z$.
  \item\label{item:lifting_3}  The obstruction to lifting the structure group of the principal
  bundle $\calg_0\xrightarrow{p} \upi_{0}(\calg_{0})\cong G$ from $\ker(p)$ to the central extension
  $\left.\calg_{1}\right|_{\{e\}\times\ker(p)}$ from \eqref{eqn:lifting_1} defines a class in $\check{H}^{2}(G,\underline{U(1)})$ which
  represents the generator under the isomorphism $\check{H}^{2}(G,\underline{U(1)})\cong \Z$.          
  \item\label{item:lifting_4} For the fibration $|\calg| \to G$ (cf.\ Proposition
        \ref{prop:separable}) the connecting homomorphism in the long
        exact homotopy sequence
        $\Z =  \pi_3(G) \to \pi_2(K(\Z,2)) = \Z $ is an isomorphism.
 \end{enumerate}
\end{lemma}
\begin{proof}
It is clear that \ref{item:lifting_1} and \ref{item:lifting_4} are equivalent since we know that $|\calg|\to G$ is fibration and $|\calg|$ has the correct homotopy in all degrees except for the third. The equivalence between \ref{item:lifting_2} and \ref{item:lifting_3} follows from the fact that the bundle gerbe in \ref{item:lifting_2} can be rewritten as a lifting bundle gerbe and then the equivalence between \ref{item:lifting_2} and \ref{item:lifting_3} is true by the fact that the class of lifting bundle gerbes is exactly the obstruction for the lifting problem. It remains to show the equivalence between \ref{item:lifting_1} and \ref{item:lifting_2}

We denote the bundle gerbe described in \ref{item:lifting_2} by $\mathfrak{B}$. By the fact that bundle gerbes are up to stable isomorphism classified by $H^3(G,\Z) = \Z$ it follows that we have $\mathfrak{B} \cong \mathcal{I}^k$ where $\mathcal{I}$ is any model for the basic bundle gerbe over $G$ and $k \in \Z$. We show that $\calg$ is a 2-group model for the string group iff $k=1$ or $-1$. Therefore it suffices to show that $||\mathcal{I}^k||$ is a 3-connected cover for $k=1,-1$ since the geometric realization is invariant under stable isomorphism of gerbes, which follows from \cite[Theorem 6.3 and Theorem 8.2]{noohi}. This last assertion now is shown in \cite[Theorem 1.2.]{BaezCransStevensonSchreiber07From-loop-groups-to-2-groups} for the choice of $\mathcal{I}$ obtained from their 2-group model. 
\end{proof}

\begin{remark}
 Considering $\String$ as a category with only identity morphisms we
 obtain a 2-group as in Example \ref{ex:bu1}. However, in this case
 $\upi_{1}\String$ is trivial. So it is not a 2-group model as defined
 above, although its geometric realization is a topological group model.
\end{remark}

\section{The string group as a 2-group}
\label{sect:the_string_group_as_a_2_group}

The previous remark shows that Lie 2-group models have more structure than topological or Lie group models for the string group. In this section we promote our Lie group model from Section \ref{sec:group} to such a Lie 2-group model. Therefore the setting will be as in Section \ref{sec:group}: $G$ is a compact simple, 1-connected Lie group and $P \to G$ is a smooth $PU(\calh)$ bundle that represents the generator $1 \in H^3(G,\Z)\cong\Z$.\\

\begin{tabsection}
 Clearly we have the central extension
 $U(1) \to U(\calh) \to PU(\calh)$. Furthermore $PU(\calh)$ acts by
 conjugation on $U(\calh)$. Using these maps we obtain a sequence
 \begin{equation}\label{bigextension}
  C^\infty (G,U(1)) \to C^\infty (P,U(\calh))^{PU(\calh)} \to \Gau(P),
 \end{equation}
 which is a central extension of Fr\'echet--Lie groups by Proposition
 \ref{prop:local_sections}.

 For the next proposition note that each smooth function
 $f \in C^\infty(G,U(1))$ is a quotient of a smooth function
 $\hat f \in C^\infty(G,\R)$ by the fact that $G$ is 1-connected. If
 we identify $U(1)$ with $\R/\Z$ we may thus identify
 $C^{\infty}(G,U(1))$ with $C^{\infty}(G,\R)/\Z$.
\end{tabsection}

\begin{lemma}\label{lem:ig}
 If $\mu$ is the Haar measure on $G$, then the map
 \begin{equation*}
  I_{G}\from  C^{\infty}(G,U(1))\to U(1),\qquad 
  I_G\left[\hat f\right] := \left[\int_G \hat f~ d\mu\right]
 \end{equation*}
 is a smooth group homomorphism. This map $I_G$ is invariant under the
 right action of $G$ on $C^\infty(G,U(1))$ which is given by left
 multiplication in the argument.
\end{lemma}
\begin{proof}
 We denote by $dI_{G} : C^\infty(G,\R)\to \R$ the map on Lie algebras that is given by
 $dI_G(f) := \int_G f~ d\mu$. First note that $dI_G$ is linear and continuous in the topology of
 uniform convergence since we have
 $|\int_{G} f~ d \mu|\leq \int_{G} |f|~ d \mu$. It thus is also
 continuous in the finer $C^{\infty}$-topology and in particular smooth.
 Furthermore it is invariant under left multiplication with $G$.
 Moreover, $dI_{G}$ factors since it maps $\Z\subset C^{\infty}(G,\R)$ to
 $\Z\subset \R$.
\end{proof}

Now we can use the group homomorphism $I_G$ to turn the smooth extension \eqref{bigextension} into a $U(1)$ extension:

\begin{definition}\label{def:definition_of_gaut}
We define
\begin{equation*}
\Gaut(P) := C^\infty (P,U(\calh))^{PU(\calh)} \times U(1) ~ \big/ \sim,
\end{equation*}
where we identify $(\varphi \cdot \mu, \lambda) \sim (\varphi, I_G(\mu)\cdot \lambda)$ for  $\mu \in C^\infty(G,U(1))$.
\end{definition}

\begin{proposition}\label{prop:contractible}
 The sequence
 \begin{equation}\label{gau_extension}
  U(1) \to \Gaut(P) \to \Gau(P)
 \end{equation}
 is a central extension of metrizable Fr{\'e}chet Lie groups and the
 space $\Gaut(P)$ is  contractible.
\end{proposition}
\begin{proof}
 By definition of $\Gaut(P)$ it is just the association of the bundle
 $C^\infty (P,U(\calh))^{PU(\calh)} \to \Gau(P)$ along the homomorphism
 $I_G: C^\infty (G,U(1)) \to U(1)$. Hence it is a smooth manifold and a
 central extension of $\Gau(P)$. More precisely we may take a locally
 smooth $C^{\infty}(G,U(1))$-valued cocycle describing the central
 extension \eqref{bigextension}. Composing this with $I_{G}$ yields then
 a locally smooth cocycle representing the central extension
 \eqref{gau_extension} (cf.\ \cite[Proposition 4.2]{Neeb02}). Since the
 modeling space is the product of the modeling space of the fiber and
 the base it is in particular Fr\'echet. In addition, $\Gaut(P)$ is
 metrizable by Lemma \ref{lem:gau_is_paracompact} and Lemma
 \ref{lem:metrizable_fiber_product}.

 Now we come to the second part of the claim. In order to show that
 $\Gaut(P)$ is weakly contractible we first define another space
 $\Gautt(P)$ using the homomorphism $ev: C^\infty(G,U(1)) \to U(1)$
 instead of $I_G$. More precisely,
 \begin{equation*}
  \Gautt(P) := C^\infty (P,U(\calh))^{PU(\calh)} \times U(1) ~ \big/ \sim_{ev}
 \end{equation*}
 where we identify
 $(\varphi \cdot \mu, \lambda) \sim_{ev} (\varphi, \mu(1)\cdot \lambda)$
 for $\mu \in C^\infty(G,U(1))$. Note that $ev$ is smooth since
 arbitrary point evaluations are so. Thus $\Gautt(P)$ is a $U(1)$
 central extension of $\Gau(P)$ and also metrizable by Lemma
 \ref{lem:metrizable_fiber_product}.

 We claim that the $\Gaut(P)$ and $\Gautt(P)$ are
 homeomorphic as spaces (not as groups). Therefore we first show that the homomorphisms $ev$ and
 $I_G$ are homotopic as group homomorphisms, i.e. there is a homotopy
 \begin{equation*}
  H: C^\infty(G,U(1)) \times [0,1] \to U(1)
 \end{equation*}
 such that each $H_t := H(-,t)$ is a Lie group homomorphism, $H_0 = ev$
 and $H_1 = I_G$. We first define the smooth map
 \begin{equation*}
  dH\from C^{\infty}(G,\R)\times[0,1]\to \R,\quad(f,t)\mapsto
  t\cdot f(1)+  (1-t)\cdot\int_{G}f~ d\mu
 \end{equation*}
 Since each $dH_{t}$ maps $\Z$ into $\Z$ it in particular induces a
 smooth group homomorphism $H_{t}$ via the identification
 $C^{\infty}(G,U(1))\cong C^{\infty}(G,R)/\Z$. Now we can use $H_t$ to
 define a $U(1)$-bundle $E$ over $\Gau(P) \times [0,1]$ by
 \begin{equation*}
  E := C^\infty (P,U(\calh))^{PU(\calh)} \times U(1) \times [0,1] ~ \big/ \sim_H
 \end{equation*}
 where we identify
 $(\varphi \cdot \mu, \lambda, t) \sim_H (\varphi, H(\mu,t) \cdot \lambda, t)$.
 Obviously $E\big|_{\Gau(P) \times 0} \cong \Gautt(P)$ and
 $E\big|_{\Gau(P) \times 1} \cong \Gaut(P)$. Thus $\Gautt(P)$ and
 $\Gaut(P)$ are isomorphic as continuous bundles \cite[Theorem 14.3.2]{Dieck08}.

 Since we now know that $\Gaut(P) \cong \Gautt(P)$, it is sufficient to
 show that $\Gautt$ is contractible. To this end we first pick a point
 $p \in P$ in the fiber over $1 \in G$. Evaluation at $p$ yields a group
 homomorphism
 \begin{equation*}
  ev: \Gau(P)=C^\infty(P,PU(\calh))^{PU(\calh)} \to PU(\calh).
 \end{equation*}
 which is a weak homotopy equivalence by \cite[Lemma 5.6]{Stolz} and
 Proposition \ref{prop:gau_homotopy_equivalence}. We now define another
 Lie group homomorphism $\Phi: \Gautt(P) \to U(\calh)$ by
 $\Phi([\varphi,\lambda]) :=  \lambda \cdot \varphi(p) $. By definition
 of $\Gautt(P)$ this is well defined and the diagram
 \begin{equation*}
  \xymatrix{
  U(1) \ar[r]\ar@{=}[d] & \Gautt(P) \ar[r]\ar[d]^\Phi & \Gau(P) \ar[d]^{ev}\\
  U(1) \ar[r] & U(\calh) \ar[r] & PU(\calh) }.
 \end{equation*}
 commutes. Since $ev$ is a weak homotopy equivalence it follows from the
 long exact homotopy sequence and the Five Lemma that also $\Phi$ is a
 weak homotopy equivalence. Therefore the weak contractibility of
 $\Gautt(P)$ is implied by the weak contractibility  of $U(\calh)$. This
 also implies contractibility of $\Gautt(P)$ by Theorem \ref{thm:CW-type}.
\end{proof}

Combining the two sequences \eqref{eqn:sting_extension} and \eqref{gau_extension} we obtain an exact sequence
\begin{equation}\label{foursequence}
1 \to U(1) \to \Gaut(P) \xrightarrow{\partial} \String \to G \to 1
\end{equation}
of Fr{\'e}chet Lie groups, where $\partial$ is the composition $\Gaut(P) \to \Gau(P) \to \String$. We furthermore define a smooth right action of $\String$ on $\Gaut(P)$ by: 
\begin{equation}\label{eqn:string_action}
[\varphi,\lambda]^f := [\varphi \circ f,\lambda] \qquad \text{ for }f \in \String \subset \Aut(P).
\end{equation}

\begin{proposition}\label{prop:crossed_module_for_string_2-group}
The action is well defined. Together with the morphism $\partial: \Gaut(P) \to \String$ this forms a smooth crossed module.
\end{proposition}

\begin{proof}
The action is well-defined since for $\varphi \in C^\infty (P,U(\calh))^{PU(\calh)}$, $\mu \in C^\infty(G,U(1))$ and $f \in \String$
we have
\begin{equation*}
[(\varphi\cdot \mu) \circ f,\lambda] = [(\varphi \circ f) \cdot (\mu \circ Q(f)),\lambda] = [\varphi \circ f,I_G(\mu \circ Q(f))\cdot \lambda] = [\varphi \op{\circ} f, I_G(\mu)\cdot\lambda]
\end{equation*}
where the last equality holds by the fact that $I_G$ is invariant under left multiplication as shown in Lemma \ref{lem:ig}.

The action of $\Aut(P)$ on $\Gau(P)\cong C^{\infty}(P,PU(\calh))^{PU(\calh)}$, given by
$\varphi^{f}:=\varphi \op{\circ} f$ is the conjugation action of $\Gau(P)$ on
itself \cite[Remark 2.8]{wockel1}. This shows that $\partial$ is equivariant and that \eqref{foursequence} and 
\eqref{eqn:string_action} define indeed a crossed module. It thus remains to show
 that the action map $\Gaut(P) \times \String \to \Gaut(P)$ is smooth. Since
$\String$ acts by diffeomorphisms it suffices to show that the restriction of the
 action map $U \times \Gaut(P)\to \Gaut(P)$ for $U$ some neighborhood of the identity in
 $\String$ is smooth. By Theorem \ref{mainthm1} we find some $U$ which is
 diffeomorphic to $\Gau(P) \times O$ for some open $O\subset G$
 with $1_{G}\in O$.
 Writing out the induced map $\Gaut(P) \times \Gau(P)\times O \to\Gaut(P)$ in local
 coordinates one sees that the smoothness of this map is implied from the
 smoothness of the action of $\Gau(P)$ on $C^{\infty}(P,U(\calh))^{PU(\calh)}$
 and the smoothness of the natural action 
$C^{\infty}(G,U(H)) \times \Diff(G) \to C^{\infty}(G,U(H))$, 
$(\varphi,f)\mapsto\varphi \op{\circ} f$
 \cite{GlocknerNeeb08Infinite-dimensional-Lie-groups-I}.
\end{proof}

\begin{definition}
Let $G$ be a compact, simple and 1-connected Lie group. Then we define $\STRING$ to be the metrizable Fr{\'e}chet Lie 2-group associated to the crossed module $\big(\Gaut(P) \xrightarrow{\partial} \String\big)$ according to example \ref{ex:bu1}. 
\end{definition}
In more detail we have 
\begin{equation*}
\big(\STRING\big)_0 := \String \qquad \text{ and } \qquad \big(\STRING\big)_1 := \Gaut(P) \rtimes \String
\end{equation*}
with structure maps given by
\begin{equation*}
s(g,f) = f \qquad t(g,f) = \partial(g)h \qquad i(f) = (1,f) \quad \text{and}\quad (g,f)\circ (g',f') = (gg',f).
\end{equation*}
From the sequence \eqref{foursequence} we obtain isomorphisms
\begin{equation}\label{isos}
\upi_0\STRING = \text{coker}(\partial) \iso G \quad \text{and} \quad \upi_1\STRING = \ker(\partial) \iso U(1).
\end{equation}
Moreover we can consider the Lie group $\String$ from Definition \ref{def:string} also as a 2-group which has only identity morphisms, see Example \ref{ex:bu1}. Then there is clearly an inclusion $\String \to \STRING$ of 2-groups.

\begin{theorem}\label{mainthm2}
The 2-group $\STRING$ together with the isomorphisms \eqref{isos}
is a smooth 2-group model for the string group (in the sense of Definition \ref{def:twomodel}). The inclusion $\String \to \STRING$ induces a homotopy equivalence 
\begin{equation*}
\String \to |\STRING|
\end{equation*}
\end{theorem}
\begin{proof}
We first want to show that the map $\String = |\String| \to |\STRING|$ is a homotopy equivalence. Therefore note that the inclusion functor $\String \to \STRING$ is given by the identity on the level of objects and by the canonical inclusion 
\begin{equation*}
\String \to \Gaut \rtimes \String
\end{equation*}
on the level of morphisms. Both of these maps are homotopy equivalences, the first for trivial reasons and the second by the fact that $\Gaut$ is  contractible as shown in Proposition \ref{prop:contractible}. Since, furthermore, both Lie-2-groups are metrizable we can apply Proposition \ref{prop:equivreal} and conclude that the geometric realization of the functor is a homotopy equivalence.

It only remains to show that $|\STRING| \to G$ is a 3-connected cover. The homotopy equivalence $\String \simeq |\STRING|$ clearly commutes with the projection to $G$. Thus the claim is a consequence of the fact that $\String$ is a smooth String group model (in particular a 3-connected cover) as shown in Theorem \ref{mainthm1}.
\end{proof}

\begin{remark}\label{bcString}
 From Remark \ref{rem:generality1} we
 obtain a crossed module $\widetilde{\Gau(P_{e}G)}\to \BCString$,
 where $\BCString$ is the restriction of the Lie group extension
 \begin{equation}
  \Gau(P)\to \Aut(P)_{0}\to \Diff(G)_{0}
 \end{equation}
 from \cite[Theorem 2.14]{wockel1} to $G\subset \Diff(G)_{0}$ and
 $\BCString \subset \Aut(P_{e}G)$ acts canonically
 $\widetilde{\Gau(P_{e}G)}:=C^{\infty}(P_{e}G,\widehat{\Omega G})^{\Omega G}$.
 As in Definition \ref{def:definition_of_gaut} we then define
 $\widehat{\Gau}(P_{e}G)$ to be associated to $\widetilde{\Gau(P_{e}G)}$
 along the homomorphism $I_{G}$. This furnishes another crossed module
 \begin{equation*}
  \widehat{\Gau}(P_{e}G) \to \BCString,
 \end{equation*}
 where the action of $\BCString \subset \Aut(P_{e}G)$ is defined in the
 same way as in as in \eqref{eqn:string_action}.
\end{remark}

\section{Comparison of string structures}
\label{sect:string_structures_and_comparison}

One reason for the importance of Lie 2-groups is that they allow for a bundle
theory analogous to bundles for Lie groups. These 2-bundles play for example a
role in mathematical physics. In particular in supersymmetric sigma models, which
are used to describe fermionic string theories, they serve as target space
background data \cite{FreedMoore,waldorf8, bunke1}. For a precise definition of 2-bundles we refer the reader to
\cite{NikolausWaldorf11} or \cite{Wockel09}.
We mainly need the following facts about smooth 2-bundles here
\begin{enumerate}
 \item For a Lie 2-group $\calg$ and a finite dimensional manifold $M$ all
       2-bundles form a bicategory $\zwoabun \calg M$ \cite[Definition
       6.1.5]{NikolausWaldorf11}.
 \item For a smoothly separable, metrizable Lie 2-group $\calg$ isomorphism
       classes of $\calg$-2-bundles are in bijection with non-abelian
       cohomology $\textnormal{\v{H}}^1(M,\calg)$ (see \cite[Definition
       2.12]{Wockel09} and \cite[Definition 3.4]{NikolausWaldorf11}, as well as
       references therein, for the definition of
       $\textnormal{\v{H}}^1(M,\calg)$) and with isomorphism classes of
       continuous $|\calg|$-bundles \cite[Theorem 4.6, Theorem 5.3.2 and
       Theorem 7.1]{NikolausWaldorf11}.
 \item For a Lie group $G$ considered as a Lie 2-group (as in example
       \ref{ex:bu1}) the definition of 2-bundles reduces to that of 1-bundles.
       More precisely we have an equivalence of bicategories
       $\bun G M \to \zwoabun G M$ where $\bun G M$ is considered as a
       bicategory with only identity 2-morphisms \cite[Example
       5.1.8]{NikolausWaldorf11}. Moreover non-abelian cohomology
       $\textnormal{\v{H}}^1(M,G)$ reduces in this case to ordinary
       \v{C}ech cohomology.
 \item\label{fact4} For a morphism of $\calg \to \calg'$ of Lie 2-groups we
       have an induced functor $\zwoabun \calg M \to \zwoabun {\calg'} M$ and
       an induced morphism
       $\textnormal{\v{H}}^1(M,\calg) \to \textnormal{\v{H}}^1(M,\calg')$. For
       a smooth weak equivalence between metrizable, smoothly separable
       2-groups the induced functor is an equivalence of bicategories.
       \cite[Theorem 6.2.2]{NikolausWaldorf11}.
\end{enumerate}

\begin{proposition}\label{prop:comparison_1-string_vs_2-string} %
The inclusion $\String \to \STRING$ induces a functor 
\begin{equation*}
\bun \String M \to \zwoabun \STRING M
\end{equation*}
which on isomorphism classes is given by the induced map
\begin{equation*}
\textnormal{\v{H}}^1\big(M,\String\big) \to \textnormal{\v{H}}^1\big(M,\STRING\big)
\end{equation*}
for each finite dimensional manifold $M$. This induced map is a bijection. 
\end{proposition}
\begin{proof}
This follows essentially from the fact that the geometric realization of the functor $\String \to |\STRING|$ is a homotopy equivalence as shown in Theorem \ref{mainthm2}. Then one knows that the induced map between isomorphism classes of continuous $\String$-bundles and $|\STRING|$-bundles is an isomorphism. Then the claim follows by the facts given above.
\end{proof}

The importance of the last proposition is that it allows us to directly compare $\String$-structures and $\STRING$-structures. We mainly built the 2-group model $\STRING$ in order to have such a comparison available. Now one can use the $\STRING$ 2-group and compare it in the world of Lie 2-groups to other smooth 2-group models and so obtain an overall comparison. We will make precise what this means in detail:
\begin{definition}
A morphism between 2-group models $\calg$ and $\calg'$ is a smooth homomorphism $f: \calg \to \calg'$ of 2-groups such that the diagrams
\begin{equation*}
\xymatrix{ 
\upi_0\calg \ar[rr]^{\upi_0 f}\ar[rd]_\sim && \upi_0\calg' \ar[ld]^\sim \\
& G &} \qquad \text{and}\qquad
\xymatrix{ 
\upi_1\calg \ar[rr]^{\upi_1 f}\ar[rd]_\sim && \upi_1\calg' \ar[ld]^\sim \\
& U(1) &}
\end{equation*}
commute.
\end{definition}
\begin{proposition}
Let $f: \calg \to \calg'$ be a morphism between metrizable, smoothly separable smooth 2-group models.
\begin{enumerate}
\item Then $f$ is automatically a smooth weak equivalence of 2-groups.
\item The geometric realization $|f|: |\calg| \to |\calg'|$ is a homotopy equivalence of topological groups. Furthermore it commutes with the projection to $G$ and the inclusion of $|\calb U(1)|$ (see Proposition \ref{prop:separable}).
\item For a manifold $M$ the induced functor
\begin{equation*}
f_*: \zwoabun{\calg} M \to \zwoabun{\calg'} M.
\end{equation*}
is an equivalence of bicategories.
\end{enumerate}
\end{proposition}
\begin{proof}
The first assertion follows from the characterization of weak equivalences given in Proposition \ref{prop:separable} and the second from Proposition \ref{morita}. The last statement is then implied by fact \ref{fact4} mentioned above.
\end{proof}

This shows that from such a morphism between 2-group models we can directly derive
comparisons between the corresponding bundle theories. Of course one should allow spans of such
morphisms. An interesting thing would be to give directly such a span connecting
our model $\STRING$ to the model given in
\cite{BaezCransStevensonSchreiber07From-loop-groups-to-2-groups}. We now give an abstract 
argument for its existence. With this it is possible to directly compare the different notions
of string bundles between any kind of 2-group model and, moreover, by Proposition
\ref{prop:comparison_1-string_vs_2-string} to the classical notions of string structures.

 The existence result will rely on the uniqueness result of Schommer-Pries from
 \cite{Schommer-Pries10Central-Extensions-of-Smooth-2-Groups-and-a-Finite-Dimensional-String-2-Group}, 
 where he considers central extensions
 \begin{equation*}
  [\calb U(1)]^{\op{fin}}\to [\Gamma]^{\op{fin}} \to G
 \end{equation*}
 of finite-dimensional smooth group stacks as models for the string group (see Remark \ref{rem:stacks}
 for the notation). A big amount in the proof of the following theorem will be about 
 making the results from 
 \cite{Schommer-Pries10Central-Extensions-of-Smooth-2-Groups-and-a-Finite-Dimensional-String-2-Group}
 also accessible in the infinite-dimensional setting (cf.\ Remark \ref{rem:stacks}).

\begin{lemma}
 If $\calg$ is a 2-group model for the string group in the sense of
 Definition \ref{def:twomodel}, then the stack $[\calg]^{\infty}$ is equivalent to one given by a finite-dimensional
 Lie groupoid $[\Gamma]^{\infty}$. Moreover, the morphisms $\calb U(1)\to \calg$ and $\calg\to G$
 induce a central extension
 \begin{equation*}
  [\calb U(1)]^{\op{fin}}\to [\Gamma]^{\op{fin}} \to G
 \end{equation*}
 of finite-dimensional smooth group stacks (and thus in the sense of \cite{Schommer-Pries10Central-Extensions-of-Smooth-2-Groups-and-a-Finite-Dimensional-String-2-Group}).
\end{lemma}

\begin{proof}
 Since $\calg_{0}\to G$ is a submersion (we will identify $\upi_{0}\calg$ with
 $G$ and $\upi_{1}\calg$ with $U(1)$ throughout) there exists a good open cover
 $\mathcal{U}:=(U_{i})_{i\in I}$ of $G$ and smooth sections
 $\sigma_{i}\from U_{i}\to \calg_{0}$. Moreover, there exist lifts
 $\widetilde{\sigma}_{ij}\from U_{ij}\to \calg_{1}$ of
 $\sigma_{i}\times \sigma_{j}\from U_{ij}\to \calg_{0}\times_{G}\calg_{0}$ and
 the maps
 $\lambda_{ijk}:= \widetilde{\sigma}_{ij}\cdot \widetilde{\sigma}_{jk}\cdot \widetilde{\sigma}_{ij}^{-1}$
 constitute a \v{C}ech 2-cocycle on $\mathcal{U}$. If
 $\Gamma_{\lambda}:= U(1)\times_{\lambda}\check{C}\mathcal{U}$ denotes the
 \v{C}ech groupoid of $\mathcal{U}$, twisted by the \v{C}ech cocycle $\lambda$
 (cf.\ \cite[Remark 2.16.]{Wockel09}), then
 \begin{equation*}
  U_{i}\ni x_{i}\mapsto \sigma_{i}(x)\in \calg_{0}\quad\text{ and }\quad
  U_{ij}\ni x_{ij}\mapsto \widetilde{\sigma}_{ij}(x)\in\calg_{1}
 \end{equation*}
 constitute a morphism $\Gamma_{\lambda}\to \calg$ of Lie groupoids. Moreover,
 a direct verification shows that it is a weak equivalence of Lie groupoids.
 Thus $[\calg]^{\infty}$ is equivalent to $[\Gamma_{\lambda}]^{\infty}$. In
 particular, all morphisms defined on $\calg$ (like the projection $\calg\to G$
 or the multiplication $\calg\times\calg\to\calg$) define respective morphisms
 of $[\Gamma_{\lambda}]^{\infty}$. Since the embedding
 $\cat{sSt}^{\op{fin}}\hookrightarrow\cat{sSt}^{\infty}$ is fully faithful (cf.\
 Remark \ref{rem:stacks}), this then
 gives to a central extension
 \begin{equation*}
  [\calb U(1)]^{\op{fin}}\rightarrow[\Gamma_{\lambda}]^{\op{fin}}\rightarrow G
 \end{equation*}
 of finite-dimensional smooth group stacks.
\end{proof}

The following fact is the same uniqueness assertion as for the finite-dimensional (but non-strict)
2-group models from
\cite{Schommer-Pries10Central-Extensions-of-Smooth-2-Groups-and-a-Finite-Dimensional-String-2-Group}.

\begin{theorem}\label{thm:comparison}
 If $\calg$ and $\calg'$ are smooth 2-group models for the string group,
 then there exists another smooth 2-group model $\calh$ and a span of
 morphisms
 \begin{equation*}
 \calg \xleftarrow{} \calh \xrightarrow{} \calg'
 \end{equation*}
 of smooth 2-group models.
\end{theorem}

\begin{proof}
 We adopt the notation of the previous lemma and of Remark \ref{rem:stacks}. Both models give rise to
 central extensions
 \begin{equation}\label{eqn:centr_extn}
  [\calb U(1)]^{\op{fin}}\rightarrow[\Gamma_{\lambda}]^{\op{fin}}\rightarrow G \quad\text{ and }\quad
  [\calb U(1)]^{\op{fin}}\rightarrow[\Gamma_{\lambda'}]^{\op{fin}}\rightarrow G 
 \end{equation}
 of smooth group stacks. The equivalence classes of such central
 extension are by \cite[Theorem
 99]{Schommer-Pries10Central-Extensions-of-Smooth-2-Groups-and-a-Finite-Dimensional-String-2-Group}
 in natural bijection with the third Segal-Mitchison topological group
 cohomology $H^{3}_{\op{SM}}(G,U(1))$.

 Now there is a spectral sequence with
 \begin{equation*}
  E_{1}^{p,q}=H^{q}(G^{p},\underline{U(1)})\Rightarrow H^{p+q}_{\op{SM}}(G,U(1))
 \end{equation*}
 (cf.\ \cite{Deligne74Theorie-de-Hodge.-III},
 \cite{WagemannWockel11A-Cocycle-Model-for-Topological-and-Lie-Group-Cohomology}).
 Since $H^{0}(G^{p},\underline{U(1)})=C^{\infty}(G^{p},U(1))$ and
 $H^{1}(G^{p},\underline{U(1)})=0$ we have that the $E_{1}$-term looks
 like
 \begin{equation*}
  \xymatrix@R=2em@C=2em{		&	0	&	\bullet	&	\bullet	&	\bullet &\bullet \\%
 	&  0\ar[r]^(0.3){d_{\op{gp}}}	&	H^{2}(G,\underline{U(1)})\ar[r]^(0.43){d_{\op{gp}}}	&	H^{2}(G\times G,\underline{U(1)})\ar[r]^(0.7){d_{\op{gp}}}	&	\cdots &\bullet \\%
  	& 0												
  & 0 	&	0	&	0 & 0 \\%
  	&	\bullet	&	\bullet
  & \cdots	\ar[r]^(0.3){d_{\op{gp}}} &	C^{\infty}(G^{3},U(1))\ar[r]^(0.7){d_{\op{gp}}} & \cdots\\%
  &				&			&			&	&
  \ar "5,1" -/d:a(-50) 2em/; "5,6" -/d:a(-50) 2em/
  \ar "5,1" -/d:a(-50) 2em/ ; "1,1" -/d:a(-50) 2em/
  }
 \end{equation*}
 Since the cohomology of the complex of globally smooth cochains (i.e.
 the cohomology of the bottom row) vanishes for compact Lie groups
 this shows that
 \begin{equation*}
  E_{\infty}^{1,2}\cong \ker\big(\xymatrix{H^{2}(G,\underline{U(1)})\ar[r]^(0.43){d_{\op{gp}}}	&	H^{2}(G\times G,\underline{U(1)})}\big)
 \end{equation*}
 and that $E_{\infty}^{p,3-p}$ vanishes if $p\neq 1$. Thus the only term of $E_{\infty}$
 that contributes to $H^{3}_{\op{SM}}(G,\underline{U(1)})$ is
 $E_{\infty}^{1,2}\cong H^{3}_{\op{SM}}(G,\underline{U(1)})$ and we
 obtain an exact sequence
 \begin{equation*}
  0\to H^{3}_{\op{SM}}(G,U(1))\to H^{2}(G,\underline{U(1)})\to \cdots .
 \end{equation*}
 On the extensions from \eqref{eqn:centr_extn} the map
 $H^{3}_{\op{SM}}(G,U(1))\to H^{2}(G,\underline{U(1)})$ is given by
 $[\Gamma_{\lambda}]^{\op{fin}}\mapsto [\lambda]$. Since this map is injective and
 $[\lambda]=[\lambda']$ by Lemma \ref{lem:lifting}, there has to be a
 morphism
 $\Phi\from[\Gamma_{\lambda}]^{\op{fin}}\to[\Gamma_{\lambda'}]^{\op{fin}}$
 of smooth group stacks that makes the diagram
 \begin{equation}\label{eqn:stack_morphism}
  \xymatrix{&[\Gamma_{\lambda}]^{\op{fin}}\ar[dr]\ar[dd]^{\Phi}\\[\calb U(1)]^{\op{fin}}\ar[ur]\ar[dr]&& ~\phantom{[/}G\hphantom{(1)]}\\&[\Gamma_{\lambda}']^{\op{fin}}\ar[ur]}
 \end{equation}
 commute. Up to here the argument is the same as is \cite[Theorem
 100]{Schommer-Pries10Central-Extensions-of-Smooth-2-Groups-and-a-Finite-Dimensional-String-2-Group}.

 By Remark \ref{rem:stacks} $\Phi$ induces also a morphism of group
 stacks
 $\Phi\from [\Gamma_{\lambda}]^{\infty}\to[\Gamma'_{\lambda}]^{\infty}$
 and thus also one $\Phi\from [\calg]^{\infty}\to[\calg']^{\infty}$, to
 which we apply the results from
 \cite{AldrovandiNoohi09Butterflies.-I.-Morphisms-of-2-group-stacks}\footnote{The
 fact that Schommer-Pries does consider a slightly different concept of
 group stacks is not relevant for this, since the morphisms are in both
 cases weakly monoidal morphisms.}. It follows from \cite[Theorem
 4.3.1]{AldrovandiNoohi09Butterflies.-I.-Morphisms-of-2-group-stacks}
 that there exists a butterfly
 \begin{equation*}
  \xymatrix@=1em{K\ar[dd]_{\partial}\ar[dr]^{\kappa}&& K'\ar[dd]^{\partial'}\ar[dl]_{\iota}\\ & H\ar[dr]^{\xi}\ar[dl]_{\zeta}\\ L && L' }
 \end{equation*}
 where $K\xrightarrow{\partial} L$ and $K'\xrightarrow{\partial'}L'$ are
 the crossed modules associated to $\calg$ and $\calg'$ (cf.\ Example
 \ref{ex:bu1}). This butterfly gives rise to a third crossed module
 \begin{equation*}
  K \times K'\to H, \quad (k,k')\mapsto \kappa(k)\cdot \iota(k')
 \end{equation*}
 with the action of $H$ on $K\times K'$ through
 $H\xrightarrow{\zeta\times \xi}L\times L'$ and the product action (cf.\
 \cite[Section
 5.4]{AldrovandiNoohi09Butterflies.-I.-Morphisms-of-2-group-stacks}).
 From this crossed module we obtain the Lie 2-group $\calh$ with the
 morphisms $\calh\to\calg$ and $\calh\to\calg'$ induced by $\zeta$ and
 $\xi$ and the projections onto $K$ and $K'$. As explained in
 \cite[Section
 5.4]{AldrovandiNoohi09Butterflies.-I.-Morphisms-of-2-group-stacks}, one
 of the morphism induces isomorphisms on $\upi_{0}$ and $\upi_{1}$, and
 since their composition does so (by the commutativity of
 \eqref{eqn:stack_morphism}) also the other is. With Proposition
 \ref{prop:separable} we see that both morphisms $\calh\to\calg$ and
 $\calh\to\calg'$ are weak equivalences. Moreover, $\calh_{0}$ can be
 taken to be metrizable (see the explicit description in the following
 remark) and one directly sees that $\calg_{1}$ is so. In particular, the
 geometric realization $|\calh|\to |\upi_{0}(\calh)|\cong G$ is a
 3-connected cover. If we define morphisms of \eqref{def:twomodel} to be
 induced by these isomorphisms on $\upi_{0}$ and $\upi_{1}$, then this
 turns $\calh$ into a 2-group model.
\end{proof}

\begin{remark}
 The Lie group $H$ in the previous proof can be made more explicit, it
 is defined such that the diagram
 \begin{equation*}
  \xymatrix{ [H]^{\infty}\ar[d] \ar[rr]&& [L']^{\infty}\ar[d]^{\pi_{\calg'}}\\[L]^{\infty}\ar[r]^{\pi_{\calg}}&[\calg]^{\infty}\ar[r]^{\Phi}& [\calg']^{\infty}}
 \end{equation*}
 is a weak pull-back in the category of stacks on the site of locally
 convex manifolds. This means that $[H]^{\infty}$ may be represented by a Lie
 subgroup of $\calg_{0}\times \calg '_{1}\times \calg'_{0}$ (see the
 definition of $\Psi$ on \cite[p.\
 715]{AldrovandiNoohi09Butterflies.-I.-Morphisms-of-2-group-stacks}).
 However, obtaining $H$ explicitly relies heavily on a detailed
 understanding of the morphism $\Phi$ and thus on the corresponding Lie
 groupoids $\Gamma_{\lambda}$, $\Gamma_{\lambda'}$ and on the morphism
 $[\Gamma_{\lambda}]^{\op{fin}}\to[\Gamma_{\lambda'}]^{\op{fin}}$ coming from the
 classification result in
 \cite{Schommer-Pries10Central-Extensions-of-Smooth-2-Groups-and-a-Finite-Dimensional-String-2-Group}.
 This seems to be a rather intractable route.

 In particular, it would be very good to have an explicit span of weak
 equivalences for the 2-group model from
 \cite{BaezCransStevensonSchreiber07From-loop-groups-to-2-groups} and
 our model, since the loop group model from
 \cite{BaezCransStevensonSchreiber07From-loop-groups-to-2-groups} is
 closely related to representations of the string group on the
 2-category of von Neumann algebras, bimodules and bimodule morphisms,
 while we expect that our model is more related to $C^{*}$-algebras
 (note for instance that $\String$ acts canonically on the continuous
 trace $C^{*}$-algebra of sections of the basic compact operator bundle
 over $G$). An explicit span would yield (via the usual pull-push
 construction for representations along spans) a relation between these
 worlds.
\end{remark}

\appendix

\section{Locally convex manifolds and Lie groups}
\label{app:locally_convex_manifolds}
\begin{tabsection}
 In this section we provide the necessary information to clarify the
 differential geometric background. If $X,Y$ are locally convex vector
 spaces and $U\subset X$ is open, then $f\from U\to Y$ is called
 \emph{continuously differentiable} if for each $v\in X$ the limit
 \begin{equation}
  df(x)(v):=\lim_{h\to 0}\frac{1}{h}(f(x+hv)-f(x))
 \end{equation}
 exists and the map $U\times X\to Y$, $(x,v)\mapsto df(x)(v)$ is
 continuous. It is called \emph{smooth} if the iterated derivatives
 $d^{n}f\from U\times X^{n}\to Y$ exist and are also continuous.
 Concepts like manifolds and tangent bundles carry over to this setting
 of differential calculus, in particular the notion of Lie groups and
 their associated Lie algebras
 \cite{GlocknerNeeb08Infinite-dimensional-Lie-groups-I}. Moreover,
 manifolds in this sense are in particular topological manifolds in the sense of
 \cite{pal66}.

 If $M,N$ are manifolds and $f\from M\to N$ is smooth, then we call $f$
 an \emph{immersion} if for each $m\in M$ there exist charts around $m$
 and $f(m)$ such that the corresponding coordinate representation of $f$
 is an inclusion of the modeling space of $M$ as a direct summand into
 the modeling space of $N$. Analogously, $f$ is called a \emph{submersion} if for
 each $m\in M$ the corresponding coordinate representation is a
 projection onto a direct summand (cf.\ \cite[\S II.2]{Lang99},
 \cite[Definition 4.4.8]{Hamilton82}).

 If $G$ is a Lie group, then a closed subgroup $H\subset G$ is called \emph{Lie
  subgroup} if it is also a submanifold. This is not automatically the
 case in infinite dimensions (cf.\ \cite[Exercise
 III.8.2]{Bourbaki98Lie-groups-and-Lie-algebras.-Chapters-1--3}).
 Moreover, if $H$ is a closed Lie subgroup, then it need not be immersed
 as the example of a non-complemented subspace in a Banach space shows.
\end{tabsection}

\begin{lemma}\label{lem:submersion}
 If $H\subset G$ is a closed Lie subgroup and $G/H$ carries an arbitrary Lie
 group structure such that $q\from G\to G/H$ is smooth, then the following are
 equivalent.
 \begin{enumerate}
  \item $G\to G/H$ admits smooth local sections around each point.
  \item $G\to G/H$ is a locally trivial bundle.
  \item $G\to G/H$ is a submersion.
 \end{enumerate}
 In any of these cases $H$ is an immersed Lie subgroup and $G/H$ carries
 the quotient topology.
\end{lemma}
 
\begin{proof}
 If $q$ admits the local smooth sections $\sigma\from U\to G$, then
 \begin{equation*}
  q^{-1}(U)\ni g\mapsto (q(g), \sigma(q(g))^{-1}\cdot g )\in U\times H
 \end{equation*}
 defines a local trivialization of $G\to G/H$. This shows equivalence of
 the first two statements and with this aid one sees also the equivalence
 with the last statement. From the second it follows in particular that
 $H\hookrightarrow G$ is an immersion. Since submersions are open, and
 since surjective open maps are quotient maps, the topology on $G/H$ has
 to be the quotient topology.
\end{proof}

\begin{definition}\label{defslipt}(cf.\ \cite[Definition 2.1]{Neeb07})
 A \emph{split Lie subgroup} of a Lie group is a closed subgroup that
 fulfills one of the three equivalent conditions of the preceding lemma.
\end{definition}

\begin{tabsection}
 Note that each immersed Lie subgroup of a Banach--Lie group is split by
 \cite[Proposition
 III.1.10]{Bourbaki98Lie-groups-and-Lie-algebras.-Chapters-1--3}. This
 implies in particular that each closed subgroup of a finite-dimensional
 Lie group is split by \cite[Theorem
 III.8.2]{Bourbaki98Lie-groups-and-Lie-algebras.-Chapters-1--3}. Also
 note that if $H$ is closed and normal and $G/H$ carries a Lie group
 structure such that $G\to G/H$ is smooth, then a single local smooth
 section can be moved around with the group multiplication to yield a
 local smooth section around each point.
\end{tabsection}

\begin{proposition}\label{prop:existence_of_fiber_products}
 If $X,Y,Z$ are manifolds, $f\from X\to Z$ is smooth and
 $g\from Y\to Z$ is a submersion then the fiber product $X\times_{Z}Y$
 exists in the category of smooth manifolds and the projection
 \begin{equation*}
  X\times_{Z}Y \to X
 \end{equation*}
 is a submersion. Moreover the identity is a submersion and the
 composition of submersions is again a submersion. That means
 submersions form a Grothendieck pretopology (see \cite[Definition
 5]{metzler}) on the category of smooth manifolds
\end{proposition}

\begin{proof}
 This is a slight generalization of \cite[4.4.10]{Hamilton82}. The proof
 of \cite[Proposition II.2.6]{Lang99}, showing that the first statement is a local
 one and of \cite[Proposition II.2.7]{Lang99}, showing this for a
 projection carry over literally to our more general setting. Moreover, the
 question of being a submersion is also local, so \cite[Proposition
 II.2.7]{Lang99} shows that $X\times_{Z}Y \to X$ is one.
\end{proof}

\begin{corollary}
The fibers of a submersion are submanifolds.
\end{corollary}

\begin{remark}\label{rem:stacks}
 The previous proposition allows us to consider two different kinds of smooth
 stacks. The first one is the usual one, 2-functors on the site of
 finite-dimensional manifolds (with the submersion topology) of the kind
 $M\mapsto \op{Bun}_{\Gamma}(M)$, where $\op{Bun}_{\Gamma}(M)$
 denotes the groupoid of principal $\Gamma$-bundles over $M$ for $\Gamma$ a
 finite-dimensional Lie groupoid. We call this a \emph{finite-dimensional}
 smooth stack. We abbreviate this with $[\Gamma]$ and if we want to emphasize
 that the site is the one of finite-dimensional manifolds we also write
 $[\Gamma]^{\op{fin}}$.
 
 The other one are 2-functors on the site of locally convex manifolds (with a
 cardinality bound for the local models to avoid set theoretical problems) of
 the kind $M\mapsto \op{Bun}_{\Gamma}(M)$ for $\Gamma$ an arbitrary (not
 necessarily finite-dimensional) Lie groupoid. We call this an
 \emph{infinite-dimensional} smooth stack. If we want to emphasize that the
 site is the one of locally convex manifolds we also write $[\Gamma]^{\infty}$.
 
 If $\Gamma$ is finite-dimensional, then we can restrict the functor
 $[\Gamma]^{\infty}$ to finite-dimensional manifolds to obtain the
 finite-dimensional smooth stack $[\Gamma]^{\op{fin}}$. Then morphisms between
 $[\Gamma]^{\infty}$ and $[\Gamma']^{\infty}$ are (by the usual argument)
 given by left-principal bibundles between $\Gamma$ and $\Gamma'$. Likewise,
 2-morphisms are given by morphisms of bibundles. Since a left-principal
 bibundle is finite-dimensional if $\Gamma$ is so, this shows that restriction
 to finite-dimensional manifolds induces an equivalence
 $\Hom([\Gamma]^{\infty},[\Gamma']^{\infty})\cong \Hom([\Gamma]^{\op{fin}},[\Gamma']^{\op{fin}})$
 if $\Gamma$ and $\Gamma'$ are finite-dimensional. In other words: the
 inclusion of finite dimensional smooth stacks into infinite dimensional smooth
 stacks, given by $[\Gamma]^{\op{fin}}\to[\Gamma]^{\infty}$ is fully faithful.
\end{remark}

\begin{tabsection}
 A manifold is called metrizable if the underlying topology is so. Note
 that metrizable is equivalent to paracompact and locally metrizable
 \cite[Theorem 1]{pal66}. Thus a Fr\'echet manifold is metrizable if and
 only if it is paracompact. Moreover, we have the following
\end{tabsection}

\begin{theorem}\label{thm:CW-type}
 A metrizable manifold has the homotopy type of a CW-complex. In
 particular, weak homotopy equivalences between metrizable manifolds are
 homotopy equivalences.
\end{theorem}

\begin{proof}
 By \cite[Theorem 14]{pal66} a metrizable manifold is dominated by
 CW-complex. By a theorem of Whitehead this implies that it has the homotopy type
 of a CW-complex (cf.\ \cite[Prop. A.11]{HatcherAT}).
\end{proof}

\section{A characterization of smooth weak equivalences}
\label{app:lie2groups}

\begin{tabsection}
 In this section we will exclusively be concerned with smoothly
 separable Lie 2-groups. Recall that for a smoothly separable Lie
 2-group $\calg$ we require among other things that $\upi_1\calg$ is a split Lie subgroup. Our main
 goal here is to prove part 1 of Proposition \ref{prop:separable}. This will be
 done in several steps.
\end{tabsection}

\begin{lemma}
\label{lem:stsect}
Let $\calg$ be a smoothly separable Lie 2-group.
Then the map $s\times t:\calg_1\to \calg_0 \times_{\upi_0\calg} {\calg}_0$ is a surjective submersion.
\end{lemma}
\begin{proof}
By definition the map $s\times t$ is a surjective map onto the submanifold
$\calg_0 \times_{\upi_0\calg} {\calg}_0$ of $\calg_0 \times {\calg}_0$.
It admits local sections because its kernel $\upi_1\calg$ is a split Lie subgroup. By
Lemma \ref{lem:submersion} this implies that it is a submersion. 
\end{proof}

\begin{proposition}
Let $f:\calg\to\calg'$ be a morphism of smoothly separable Lie 2-groups inducing an isomorphism on
$\upi_1$. Then $f$ is smoothly fully faithful, i.e., 
\[
\xymatrix{
{\calg}_1 \ar[r]^{f} \ar[d]_{s \times t} & {\calg'}_1 \ar[d]^{s \times t} \\
{\calg}_0 \times {\calg}_0 \ar[r]_-{f \times f} & {\calg'}_0 \times {\calg'}_0
}
\]
is a pullback diagram of Lie groups.
\end{proposition}

\begin{proof}
 It is clear that this is a pullback diagram of groups by the general
 theory of 2-groups. Let $\mathcal{H}$ be a Lie group and consider the
 diagram
 \[
  \xymatrix{ \mathcal{H} \ar@/^/[drr]^a \ar@/_/[ddr]_b \ar@{-->}[dr]^h
  && \\
  & {\calg}_1 \ar[r]^{f} \ar[d]_{s \times t} & {\calg'}_1 \ar[d]^{s
  \times t} \\
  & {\calg}_0 \times {\calg}_0 \ar[r]_-{f \times f} & {\calg'}_0 \times
  {\calg'}_0 }
 \]
 where $a,b$ are morphisms of Lie groups. We have to show that the
 unique map $h:\mathcal{H}\to\calg_1$ supplied by the pullback of groups
 is also smooth. By Lemma \ref{lem:stsect} there exists a smooth local
 section $\gamma \from U\to \calg_{1}$ of $s\times t$, defined on a
 neighborhood $U\subset \calg_{0}\times_{\upi_{0}}\calg_{0}$ of the identity.
 Since $b$ maps to $\calg_{0}\times_{\upi_{0}}\calg_{0}$, $V:=b^{-1}(U)$
 is an open neighborhood of the identity in $\calh$.

 We now observe that
 \begin{equation*}
  h'\from V\to \calg_{1},\quad  
  x\mapsto 
  \gamma(b(x))\cdot (\upi_{1}f_{1})^{-1}( f_{1}(\gamma(b(x)))^{-1} \cdot a(x))
 \end{equation*}
 is smooth since
 $f_{1}(\gamma(b(x)))^{-1} \cdot a(x)\in \upi_{1}\calg' $ and
 $f_{1}$ restricts to a diffeomorphism $\upi_{1}\calg\to\upi_{1}\calg'$.
 It satisfies $f_{1} \op{\circ} h'=\left.a\right|_{V}$, and we also
 have $(s\times t) \op{\circ} h'=b$ since $\gamma$ is a section of
 $s\times t$. Thus $h$ coincides with $h'$ on $V$, showing that $h$ is a
 smooth homomorphism of Lie groups.
\end{proof}

\begin{proposition}
 Let $f:\calg\to\calg'$ be a morphism of smoothly separable Lie 2-groups
 inducing an isomorphism on $\upi_0$. Then $f$ is smoothly essentially
 surjective, i.e., the morphism
 \[
  s\circ\mathrm{pr}_2:\calg_0 {~}_{f_0}\!\times_{t} {\calg'}_1 \to
  \calg_{0}'
 \]
 is a smooth surjective submersion.
\end{proposition}
\begin{proof}
 Surjectivity is clear because $f$ is surjective on $\upi_0$. To see
 that $s\circ\mathrm{pr}_2$ is a submersion we will construct a local
 smooth section. Since the map $p:\calg_0\to\upi_0\calg$ is a submersion
 there exists a local section $\sigma\from U\to \calg_{0}$ of $p$.
 For brevity let us denote the ``roundtrip'' map, restricted to
 $V:=p'^{-1}(\upi_{0}f(U))$ as
 $R=f_0\circ\sigma\circ(\upi_0f)^{-1}\circ p'$.
 For $x\in V$ we then have $x\cong R(x)$ and thus
 $(x,R(x))\in \calg_{0}'\times_{\upi_{0}\calg}\calg_{0}' $. Now there
 exists a local smooth section $\tau\from W\to \calg'_{1}$ of
 $s'\times t'$ for $W\subset V\times_{\upi_{0}\calg'} V$ open.
 Then
 \begin{eqnarray*}
  S:(\left.\id_{\calg_{0}'}\right|_{V}\times R)^{-1}(W) &\to & \calg_0
  {~}_{f_0}\!\times_{t} {\calg'}_1\\
  x &\mapsto & (\sigma((\upi_{0}f)^{-1}(p'(x))),\tau(x,R(x)))
 \end{eqnarray*}
 is the required section since we have
 \begin{equation*}
 f_{0}(\sigma((\upi_{0}f)^{-1}(p'(x))))=R(x)=t(\tau(x,R(x)))
 \end{equation*} 
and
 $s(\tau(x))=x$.
\end{proof}

\begin{corollary}
If $f:\calg\to\calg'$ is a morphism of smoothly separable Lie 2-groups inducing isomorphisms on
$\upi_0$ and $\upi_1$ then $f$ is a weak equivalence.
\end{corollary}

The converse of the first part of Proposition \ref{prop:separable} also holds:

\begin{proposition}
 A smooth weak equivalence $f:\calg\to\calg'$ of smoothly separable Lie
 2-groups induces isomorphisms on $\upi_0$ and $\upi_1$.
\end{proposition}
\begin{proof}
 Since $f$ is in particular an equivalence of the underlying categories
 in the set-theoretic sense, it is clear that its induced morphisms
 $\upi_0f:\upi_0\calg\to\upi_0\calg'$ and
 $\upi_1f:\upi_1\calg\to\upi_1\calg'$ are group isomorphisms. From the
 diagram
 \begin{equation}
  \label{diag:pb1}
  \xymatrix{\calg_0 {~}_{f_0}\!\times_{t} {\calg'}_1 \ar[r]^{\mathrm{pr}_2} \ar[d] & \calg'_1 \ar[d]^{s}\\
  \calg_0 \ar[r]^{f_0} \ar[d]_p & \calg_0' \ar[d]^{p'}\\
  \upi_0\calg \ar[r]_{\upi_0f} & \upi_0\calg'_0}
 \end{equation}
 we see that $\upi_0f$ is smooth since we can pick a local section
 $\sigma:\upi_0\calg\to\calg_0$ of the submersion
 $p:\calg_0\to\upi_0\calg$, which shows that locally
 \[
  \upi_0f=p'\circ f_0\circ\sigma.
 \]
 To see that $(\upi_0f)^{-1}$ is smooth as well we choose a local
 section $\sigma':\upi_0\calg'\to\calg'_0$. Since we know that
 $s\circ\mathrm{pr}_2:\calg_0 {~}_{f_0}\!\times_{t} {\calg'}_1\to\calg'_0$
 is a submersion, we can also choose a section $\tau$ for that map, and
 composing $\tau\circ\sigma'$ with the projection to $\calg_0$ and
 finally to $\upi_0\calg$ coincides with $(\upi_0f)^{-1}$ which is
 therefore smooth.

 To see that $\upi_1f$ is a diffeomorphism we use the fact that the
 diagram of part 2 of the definition of a smooth weak equivalence is a
 pullback diagram. This implies in particular that the restriction of
 $f_{1}$ to the fiber over $(1,1)$, which is the submanifold
 $\upi_{1}\calg$, is a smooth bijective map.
 That its inverse is also smooth
 follows from the universal property of the pullback: there exists a
 unique smooth map $H:\upi_1\calg'\to\upi_1\calg$ that makes the
 diagram
 \[
  \xymatrix{\upi_1\calg' \ar@/_/[ddr] \ar@{-->}[dr]^{H}
  \ar@/^/[drr]^{\id} &&\\
  & \upi_1\calg \ar[d]_{s\times t} \ar[r]^{f_1} & \upi_1\calg'
  \ar[d]^{s\times t}\\
  & (1,1) \ar[r]_{f_0\times f_0} & (1,1)}
 \]
 commute, so $f_1\circ H=\id_{\upi_1\calg'}$ which means that
 $H$ is the inverse of $f_{1}$ on $\upi_{1}\calg'$, which thus is smooth.
\end{proof}

This concludes the proof of the first part of Proposition \ref{prop:separable}.

\bibliographystyle{new}

\end{document}